\documentclass[a4paper]{article}

\usepackage{RR}
\usepackage[french]{babel}

\usepackage{epsfig}
\usepackage{amssymb}
\usepackage{graphics}
\usepackage{hyperref}

\usepackage{mathrsfs}
\usepackage[latin1]{inputenc}
\usepackage{fancyhdr}
\usepackage{textcomp}
\usepackage{enumitem}
\usepackage{lineno}

\usepackage{latexsym} 
\usepackage{amsmath,amssymb,txfonts}
\usepackage{graphicx}
\usepackage[usenames]{color}

\usepackage{wasysym}

\newcommand{\PartIntSup}[1]{\left\lceil #1\right\rceil}
\newcommand{\PartIntInf}[1]{\left\lfloor #1\right\rfloor}

\usepackage{pb-diagram}

\def\lgem{\hbox{l\kern-.08em\raise .7ex\hbox{.}\kern-.08eml}}
\def\Lgem{\hbox{L\kern-.08em\raise .7ex\hbox{.}\kern-.08emL}}

  \renewenvironment{thebibliography}[1]{%
    \begin{oldthebibliography}{#1}%
      \setlength{\parskip}{0ex}%
      \setlength{\itemsep}{2.0ex}%
  }%
  {%
    \end{oldthebibliography}%
  }

\newtheorem{lemma} {Lemma}[section]

\newtheorem{example}{Example}[section]

\def \cB {{\cal B}}

\def \cD {{\cal D}}

\def\endproofbox{\hskip 1.3em\hfill\rule{6pt}{6pt}}
\newenvironment{proof}%
{%
\noindent{\bf Proof.}
}%
{%
\quad\hfill\endproofbox\vspace*{2ex} }

\def\Ca{C}
\def\Cb{C'}

\newtheorem{thm}{Theorem}[section]
\newtheorem{cor}[thm]{Corollary}

\newtheorem{rmk}[thm]{Remark}

\def\ben{\begin{enumerate}%
\itemsep 1pt plus 1pt minus 1pt}
\def\een{\end{enumerate}}
\def\bit{\begin{itemize}%
\itemsep 1pt plus 1pt minus 1pt}
\def\eit{\end{itemize}}

\RRtitle{Groupage de trafic sur deux périodes avec coefficient 4}

\RRetitle{Drop cost and wavelength optimal two-period grooming with ratio 4}

\RRauthor{Jean-Claude Bermond \thanks{Mascotte joint Project, INRIA/CNRS-I3S/UNS, 2004 route des
Lucioles, B.P. 93 - F-06902, Sophia-Antipolis, France} \and Charles J. Colbourn
\thanks{Computer Science and Engineering, Arizona State University, Tempe, AZ, 85287-8809, U.S.A.}
\and Lucia Gionfriddo \thanks{Matematica e Informatica, Universit\`{a} di Catania, viale A. Doria,
6, 95125 Catania, Italy} \and Gaetano Quattrocchi \thanks{Matematica e Informatica, Universit\`{a}
di Catania, viale A. Doria, 6, 95125 Catania, Italy} \and Ignasi Sau
\thanks{Mascotte joint Project, INRIA/CNRS-I3S/UNS,
2004 route des Lucioles, B.P. 93 - F-06902, Sophia-Antipolis, France}
\thanks{Graph Theory and Combinatorics group, MA4, Universitat Polit\`ecnica de Catalunya, Barcelona, Spain}}

\RRdate{November 2009}

\RRabstract{We study grooming for two-period optical networks, a variation of the traffic grooming
problem for WDM ring networks introduced by Colbourn, Quattrocchi, and
  Syrotiuk. In the two-period grooming problem, during the first period of time, there is all-to-all
  uniform traffic among $n$ nodes, each request using $1/C$ of the
  bandwidth; and during the second period, there is all-to-all
  uniform traffic only among a subset $V$ of $v$ nodes, each request
now  being allowed to use $1/C'$ of the bandwidth, where $C' < C$. We determine the minimum drop
cost (minimum number of ADMs) for any $n,v$ and  $C=4$ and $C' \in \{1,2,3\}$. To do this, we use
tools of graph decompositions. Indeed the two-period grooming problem corresponds to minimizing the
total number of  vertices in a partition  of the edges of the complete graph $K_n$ into subgraphs,
where each subgraph has at most $C$ edges  and where  furthermore  it contains at most $C'$ edges
of the complete graph on $v$ specified vertices. Subject to the condition that the two-period
grooming has the least drop cost, the minimum number of wavelengths required is also determined in
each case.}

\RRresume{Nous \'etudions dans cet article une variante introduite par Colbourn, Quattrocchi et
Syrotiuk du probl\`eme de groupage de trafic dans les réseaux en anneaux WDM ; le groupage de
trafic sur deux p\'eriodes. Au cours de la première période de temps le trafic considéré est
\emph{all-to-all} uniforme entre $n$ sommets, chaque requête utilisant $1 / C $ de la bande
passante; pendant la deuxième période, le trafic est
 \emph{all-to-all} uniforme mais entre les sommets d'un sous-ensemble $V$ de taille $ v < n$, chaque requête
étant alors autorisée à utiliser $1/C'$ de la bande passante, où $ C' <C $. Nous déterminons le
coût minimum (nombre minimum de ADMs), pour tout $ n, v $ et $ C = 4 $ et $ C '\in \{1,2,3 \}$.
Pour ce faire, nous utilisons les décompositions de graphes. En effet, le probl\`eme du groupage de
trafic sur deux p\'eriodes revient à minimiser le nombre total de sommets d'une partition des
arêtes du graphe complet $K_n$ en sous-graphes, où chaque sous-graphe a au plus $C$ arêtes et
contient au plus $C '$ arêtes du graphe complet sur les $v$ sommets spécifiés. Nous d\'eterminons
aussi, pour les valeurs ci-dessus, le coût minimum d'une solution utilisant  le nombre minimal
requis de longueurs d'onde.}

\RRkeyword{traffic grooming, SONET ADM, optical networks, graph decomposition, design theory.}

\RRmotcle{R\'eseaux optiques, SONET, groupage de trafic, ADM, d\'ecomposition des graphes.}

\RRprojet{MASCOTTE} \RRtheme{\THCom} \URSophia

\begin{document}
\RRNo{7101}
\makeRR

\section{Introduction}

Traffic grooming is the generic term for packing low rate signals into higher speed streams (see
the surveys~\cite{BeCo06,DuttaR02,ModianoL01,SimmonsGS99,ZhMu03}). By using traffic grooming, one
can bypass the electronics in the nodes which are not sources or destinations of traffic, and
therefore reduce the cost of the network. Here we consider unidirectional SONET/WDM ring networks.
In that case, the routing is unique and we have to assign to each request between two nodes a
wavelength and some bandwidth on this wavelength. If the traffic is uniform and if a given
wavelength has capacity for at least $C$ requests, we can assign to each request at most
$\frac{1}{C}$ of the bandwidth.  $C$ is known as the \emph{grooming ratio} or the \emph{grooming
factor}. Furthermore if the traffic requirement is symmetric, it can be easily shown (by exchanging
wavelengths) that there always exists an optimal solution in which the same wavelength is given to
each pair of symmetric requests. Thus without loss of generality we assign to each pair of
symmetric requests, called a \emph{circle}, the same wavelength. Then each circle uses
$\frac{1}{C}$ of the bandwidth in the whole ring. If the two end-nodes of a circle are $i$ and $j$,
we need one ADM at node $i$ and one at node $j$. The main point is that if two requests have a
common end-node, they can share an ADM if they are assigned the same wavelength. For example,
suppose that we have symmetric requests between nodes $1$ and $2$, and also between $2$ and $3$. If
they are assigned two different wavelengths, then we need 4 ADMs, whereas if they are assigned the
same wavelength we need only 3 ADMs.

The so called traffic grooming problem consists in minimizing the total number of ADMs to be used,
in order to reduce the overall cost of the network.

Suppose we have a ring with $4$ nodes $\{0,1,2,3\}$ and all-to-all uniform traffic. There are
therefore 6 circles (pairs of symmetric requests) $\{i,j\}$ for $0 \leq i < j \leq 3$. If there is
no grooming we need 6 wavelengths (one per circle) and a total of 12 ADMs. If we have  a grooming
factor $C=2$, we can put on the same wavelength two circles, using 3 or 4 ADMs according to whether
they share an end-node or not. For example we can put together $\{1,2\}$ and $\{2,3\}$ on one
wavelength; $\{1,3\}$ and $\{3,4\}$ on a second  wavelength, and $\{1,4\}$ and $\{2,4\}$ on a third
one, for a total of 9 ADMs. If we allow a grooming factor $C=3$, we can use only 2 wavelengths. If
we put together on one wavelength $\{1,2\}$, $\{2,3\}$, and $\{3,4\}$ and on the other one
$\{1,3\}$, $\{2,4\}$, and $\{1,4\}$ we need 8 ADMs (solution \emph{a}); but we can do better by
putting on the first wavelength $\{1,2\}$, $\{2,3\}$ and $\{1,3\}$ and on the second one $\{1,4\}$,
$\{2,4\}$ and $\{3,4\}$, using 7 ADMs (solution \emph{b}).

Here we study the problem for a unidirectional SONET ring with $n$ nodes, grooming ratio $C$, and
all-to-all uniform unitary traffic. This problem has been modeled as a graph partition problem in
both \cite{BermondCICC03} and \cite{GHLO03}. In the all-to-all case the set of requests is modelled
by the complete graph $K_n$. To a wavelength $k$ is associated a subgraph $B_k$ in which each edge
corresponds to a pair of symmetric requests (that is, a circle) and each node to an ADM. The
grooming constraint, i.e. the fact that a wavelength can carry at most $C$ requests, corresponds to
the fact that the number of edges $|E(B_k)|$ of each subgraph $B_k$ is at most $C$. The cost
corresponds to the total number of vertices used in the subgraphs, and the
 objective is therefore to minimize this number.\\

\noindent \textsc{Traffic Grooming in the Ring} 
\begin{itemize}
\item[]\textbf{\textsc{Input:}} Two integers $n$ and $C$.
\item[]\textbf{\textsc{Output:}} Partition $E(K_n)$
into
  subgraphs $B_k$, $1\leq k \leq \Lambda$, s.t.
  $|E(B_k)|\leq C$ for all $k$.

\item[]\textbf{\textsc{Objective:}} Minimize
$\sum_{k=1}^\Lambda |V(B_k)|$.
  \end{itemize}

In the example above with $n=4$ and $C=3$,  solution \emph{a} consists of a decomposition of $K_4$
into two paths with four vertices  $[1,2,3,4]$ and $[2,4,1,3]$, while solution \emph{b} corresponds
to a decomposition into a triangle $(1,2,3)$ and a star with the edges $\{1,4\}$,  $\{2,4\}$, and
$\{3,4\}$.

With the all-to-all set of requests, optimal constructions for a given grooming ratio $C$ have been
obtained using tools of graph and design theory~\cite{handbook}, in particular for grooming ratio
$C=3$~\cite{BermondC03}, $C=4$~\cite{BermondCICC03,Hu02}, $C=5$~\cite{BermondCLY04},
$C=6$~\cite{BermondCCGLM}, $C=7$~\cite{cfgll} and $C\geq N(N-1)/6$~\cite{BCM03}.

Graph decompositions have been extensively studied for other reasons as well.  See \cite{Bryant}
for an  excellent survey, \cite{ColbournR99} for relevant material on designs with blocksize three,
and \cite{handbook} for terminology in design theory.

Most of the papers on grooming deal with a single (static) traffic matrix. Some articles consider
variable (dynamic) traffic, such as finding a solution which works for the maximum traffic
demand~\cite{BeMo00,ZZM03} or for all request graphs with a given maximum degree~\cite{MuSa08}, but
all keep a fixed grooming factor. In~\cite{cqsupper} an interesting variation of the traffic
grooming problem, grooming for two-period optical networks, has been introduced in order to capture
some dynamic nature of the traffic. Informally, in the two-period grooming problem each time period
supports different traffic requirements. During the first period of time there is all-to-all
uniform traffic among $n$ nodes, each request using $1/C$ of the bandwidth; but during the second
period there is all-to-all traffic only among a subset $V$ of $v$ nodes, each request now being
allowed to use a larger fraction of the bandwidth, namely $1/C'$ where $C' < C$.

Denote by $X$ the subset of $n$ nodes. Therefore the two-period
grooming problem can be expressed as follows:\\

\noindent \textsc{Two-Period Grooming in the Ring}

\begin{itemize}
\item[]\textbf{\textsc{Input:}} Four integers $n$, $v$, $C$, and $C'$.
\item[]\textbf{\textsc{Output:}} A partition (denoted $N(n,v;C,C')$) of $E(K_n)$ into
  subgraphs $B_k$, $1\leq k \leq \Lambda$, such that for all
  $k$,
  $|E(B_k)|\leq C$, and $|E(B_k)\cap (V\times V)|\leq C'$, with
  $V\subseteq X$, $|V|=v$.

\item[]\textbf{\textsc{Objective:}} Minimize $\sum_{k=1}^\Lambda
|V(B_k)|$.
  \end{itemize}

Following~\cite{cqsupper}, a grooming is denoted by $N(n,C)$. When the grooming $N(n,C)$ is
\emph{optimal}, i.e. minimizes the total ADM cost, then the grooming is denoted by
$\mathscr{ON}(n,C)$. Whether general or optimal, the drop cost of a grooming is denoted by $cost\
N(n,C)$ or $cost\ \mathscr{ON}(n,C)$, respectively.

A grooming of a two-period network $N(n,v;C,C')$ with grooming ratios $(C,C')$ coincides with a
graph decomposition $(X,\mathcal{B})$ of $K_n$ (using standard design theory terminology,
$\mathcal{B}$ is the set of all the \emph{blocks} of the decomposition) such that $(X,\mathcal{B})$
is a grooming $N(n,C)$ in the first time period, and $(X,\mathcal{B})$ faithfully embeds a graph
decomposition of $K_v$ such that $(V,\mathcal{D})$ is a grooming $N(v,C')$ in the second time
period. Let $V \subseteq X$. The graph decomposition $(X,\cB)$ {\em embeds} the graph decomposition
$(V,\cD)$ if there is a mapping $f: \cD \to \cB$ such that $D$ is a subgraph of $f(D)$ for every $D
\in \cD$. If $f$ is injective (i.e., one-to-one), then $(X,\cB)$ {\em faithfully embeds} $(V,\cD)$.
This concept of faithfully embedding has been explored in \cite{ColbournLQ03,Quattrocchi02}.

 We
use the notation
 $\mathscr{ON}(n, v; C, C')$ to denote an optimal
grooming $N(n, v; C, C')$.

As it turns out, an $\mathscr{ON}(n, v; C, C')$ does not always coincide with an
$\mathscr{ON}(n,C)$. Generally we have $cost\ \mathscr{ON}(n, v; C, C') \geq cost\
\mathscr{ON}(n,C)$ (see Examples \ref{ex:2} and \ref{ex:3}). Of particular interest is the case
when $cost\ \mathscr{ON}(n, v; C, C')$ = $cost\ \mathscr{ON}(n,C)$ (see Example \ref{ex:1}).

\begin{example}
\label{ex:1} Let $n=7$, $v=4$, $C=4$. Let $V=\{0,1,2,3\}$ and $W=\{a_0,a_1,a_2\}$. An optimal
decomposition is given by the three triangles $(a_0,0,1)$, $(a_1,1,2)$, and $(a_2,2,3)$, and the
three 4-cycles $(0,2,a_0,a_1)$, $(0,3,a_0,a_2)$, and $(1,3,a_1,a_2)$, giving a total cost of 21
ADMs.

This solution is valid and optimal for both $C'=1$ and $C'=2$, and it is optimal for the classical
\textsc{Traffic Grooming in the Ring} problem when $n=7$ and $C=4$. Therefore, $cost\
\mathscr{ON}(7, 4; 4, 1) = cost\ \mathscr{ON}(7, 4; 4, 2)= cost\ \mathscr{ON}(7,4) = 21$.
\end{example}

 \begin{example}
\label{ex:2} Let $n=7$, $v=5$, $C=4$, and $C'=2$. Let $V=\{0,1,2,3,4\}$ and $W=\{a_0,a_1\}$. We see
later that an optimal decomposition is given by the five kites $(a_0,1,2;0)$, $(a_0,3,4;1)$, $(a_1,
1, 3;2)$, $(a_1, 2, 4; 0)$ and $(a_0,a_1,0;1)$, plus the edge $\{0,3\}$, giving a total cost of 22
ADMs. So $cost\ \mathscr{ON}(7, 5; 4, 2) = 22$. Note that this decomposition is not a valid
solution for $C'=1$, since there are subgraphs containing more than one edge with both end-vertices
in $V$.
\end{example}

\begin{example}
\label{ex:3} Let $n=7$, $v=5$, $C=4$, and $C'=1$. Let again $V=\{0,1,2,3,4\}$ and $W=\{a_0,a_1\}$.
We see later that an optimal decomposition is given by the four $K_3$s $(a_0,1,2)$, $(a_0,3,4)$,
$(a_1,0,3)$, and $(a_1,2,4)$, the $C_4$ $(0,1,a_1,a_0)$, plus the five edges $\{0,4\}$, $\{1,3\}$,
$\{0,2\}$, $\{1,4\}$, and $\{2,3\}$, giving a total cost of 26 ADMs. So $cost\ \mathscr{ON}(7, 5;
4, 1) = 26$.
\end{example}

C.J. Colbourn, G. Quattrocchi and V.R. Syrotiuk \cite{cqsupper,cqslower} completely solved the
cases when $C=2$ and $C=3$ ($C'= 1$ or $2$). In this article we determine the minimum drop cost of
an $N(n,v;4,C')$ for all $n \geq v \geq 0$ and $C' \in \{1,2,3\}$.

We are also interested in determining the minimum number of wavelengths, or \emph{wavecost},
required in an assignment of wavelengths to a decomposition. Among the $\mathscr{ON}(n,4)$s one
having the minimum wavecost  is  denoted by $\mathscr{MON}(n,4)$, and the corresponding minimum
number of wavelengths by $wavecost \mathscr{MON}(n,4)$. We characterize the
$\mathscr{ON}(n,v;C,C')$ whose wavecost is minimum among all $\mathscr{ON}(n,v;C,C')$s  and denoted
one by $\mathscr{MON}(n,v;C,C')$; the wavecost is itself denoted by $wavecost
\mathscr{MON}(n,v;C,C')$.

We deal separately with each value of $C' \in \{1,2,3\}$.  Table \ref{costform} summarizes   the
cost formulas for $n=v+w > 4$.

\begin{table}[htbp]
\begin{itemize}

\item[] $cost\ \mathscr{ON}(v+w,v;4,1) =\left\{\begin{array}{ll}

{v+w \choose 2}&\mbox{ if $v \leq w+1$}\\
\\
{v+w \choose 2} + {v \choose 2} - \PartIntInf{\frac{vw}{2}}&\mbox{ if $v \geq w+1$}

\end{array}\right.$\\

\vspace{0.7cm}

\item[]

$cost\ \mathscr{ON}(v+w,v;4,2) =\left\{\begin{array}{ll}

{v+w \choose 2}&\mbox{ if $v \leq 2w$}\\
\\
{v+w \choose 2} + \PartIntSup{\frac{1}{2}{v \choose 2}} - \frac{vw}{2}+
\delta &\mbox{ if $v > 2w$ and $v$ even}\\
\\

\ \ \ \ \ \ \ \mbox{ where }\delta=\left\{\begin{array}{cl}

1 & \mbox{ if }w=2, \mbox{ or }\\
& \mbox{ if }w=4\mbox{ and }\\
&\ \ \ \ \ v\equiv 0 \pmod{4}\\
0 & \mbox{ otherwise}
\end{array}\right.\\

\\
{v+w \choose 2} + \PartIntSup{\frac{1}{2}\left(  {v \choose 2} - vw -\PartIntSup{\frac{w}{2}}
\right)} + \delta&\mbox{ if $v >
2w$ and $v$ odd}\\
\\

\ \ \ \ \ \ \ \mbox{ where }\delta=\left\{\begin{array}{cl}

1& \mbox{ if }w=3 \mbox{  and }\\
&\ \ \ \ \ v\equiv 3 \pmod{4}
\\
0& \mbox{ otherwise}
\end{array}\right.\\
\end{array}\right.$
\vspace{0.7cm}

\item[] $cost\ \mathscr{ON}(v+w,v;4,3)={v+w \choose 2}$

\end{itemize}

\caption{Cost formulas for $n=v+w > 4$.} \label{costform}
\end{table}

\section{Notation and Preliminaries}
\label{sec:prelim}

We establish some graph-theoretic notation to be used throughout. We denote the edge between $u$
and $v$ by $\{u,v\}$.
 $K_n$ denotes a complete
graph on $n$ vertices and  $K_X$  represents the complete graph on the vertex set $X$. A triangle
with edges $\{\{x,y\},\{x,z\},\{y,z\}\}$ is denoted by $(x,y,z)$. A 4-cycle with edges
$\{\{x,y\},\{y,z\},\{z,u\},\{u,x\}\}$ is denoted by $(x,y,z,u)$. A {\em kite} with edges
$\{\{x,y\},\{x,z\},\{y,z\},\{z,u\}\}$ is denoted by $(x,y,z;u)$. The groomings to be produced also
employ paths; the path on $k$ vertices $P_k$ is denoted by $[x_1,\dots,x_k]$ when it contains edges
$\{x_i,x_{i+1}\}$ for $1 \leq i < k$.  Now let $G=(X,E)$ be a graph.  If $|X|$ is even, a set of
$|X|/2$ disjoint edges in $E$ is a {\em 1-factor}; a partition of $E$ into 1-factors is a {\em
1-factorization}.  Similarly, if $|X|$ is odd, a set of $(|X|-1)/2$ disjoint edges in $E$ is a {\em
near 1-factor}; a partition of $E$ into near 1-factors is a {\em near 1-factorization}. We also
employ well-known results on partial triple systems and group divisible designs with block size
three; see \cite{ColbournR99} for background.

The vertices of the set $V$ are the integers modulo $v$ denoted by $0,1,\ldots,v-1$. The vertices
not in $V$, that is in $X\setminus V$, forms the set $W$ of size $w =n-v$ and is denoted by
 $a_0,\ldots,a_{w-1}$, the indices being taken
modulo $w$.

Among graphs with three or fewer edges (i.e., when $\Ca = 3$),the only graph with the minimum ratio
(number of vertices over the number of edges) is the triangle.  For $C=4$ three different such
graphs have minimum ratio 1: the triangle, the 4-cycle, and the kite.  This simplifies the problem
substantially.  Indeed, in contrast to the lower bounds in \cite{cqslower}, in this case the lower
bounds arise from easy classification of the edges on $V$. We recall the complete characterization
for optimal groomings with a grooming ratio of four:

\begin{thm}\label{lavbermcoudhu}{\rm \cite{BermondCICC03,Hu02}}
$cost\,\mathscr{ON}(4,4)=7$ and, for $n\geq 5$, $cost\,\mathscr{ON}(n,4)={n\choose 2}$. Furthermore
a $\mathscr{MON}(4,4)$ employs two wavelengths and can be realized by a kite and a $P_3$ (or a
$K_3$ and a star), and a $\mathscr{MON}(n,4)$, $n\geq 5$, employs
$\left\lceil\frac{n(n-1)}{8}\right\rceil$ wavelengths and can be realized by $t$ $K_3$s and
$\left\lceil\frac{n(n-1)}{8}-t\right\rceil$ 4-cycles or kites, where $$\ t=  \left \{
\begin{array}{ll}
  0\ \ \ \mbox{if}\ \  n\equiv 0,1\pmod{8}\\
  1\ \ \ \mbox{if}\ \  n\equiv 3,6\pmod{8}\\
  2\ \ \ \mbox{if}\ \  n\equiv 4,5\pmod{8}\\
  3\ \ \ \mbox{if}\ \  n\equiv 2,7\pmod{8}\\
\end{array}
\right..$$
\end{thm}

 In order to unify the treatment of the lower bounds, in a decomposition
$N(v+w,v;4,C')$ for $C'\in\{1,2\}$, we call an edge with both ends in $V$ {\sl neutral} if it
appears in a triangle, 4-cycle, or kite; we call it {\sl positive} otherwise.  An edge  with one
end in $V$ and one in $W$ is a  {\sl cross edge}.

\begin{lemma}\label{neutral}
 \begin{enumerate}
\item In an
$N(v+w,v;4,C')$ with $C'\in\{1,2\}$, the number of neutral edges is at most $\frac{1}{2} C' vw$.
\item When $v$ is odd and $C'=2$, the number of neutral edges is at most $ vw - \frac{w}{2}$.
\end{enumerate}
\end{lemma}
\begin{proof} Every neutral
edge appears in a subgraph having at least two cross edges. Thus the number of subgraphs containing
one or more neutral edges is at most $\frac{1}{2}vw$.  Each can contain at most $C'$ neutral edges,
and hence there are at most $\frac{1}{2} C'vw$ neutral edges.  This proves the first statement.

Suppose now  that $C'=2$ and  $v$ is odd.  Any subgraph containing two neutral edges employs
exactly  two cross edges incident to the same vertex in $W$. Thus the number $\alpha$ of such
subgraphs is at most $\frac{1}{2}w(v-1)$. Then remaining neutral edges must arise (if present) in
triangles, kites, or 4-cycles that again contain two cross edges but only one neutral edge; their
number, $\beta$, must satisfy $\beta \leq \frac{vw}{2} - \alpha$.  Therefore the number of neutral
edges, $2\alpha+\beta$, satisfies $ 2\alpha+\beta \leq \frac{1}{2}w(v-1) + \frac{vw}{2} = vw -
\frac{w}{2}$.
\end{proof}

When $C=3$ there are strong interactions among the decompositions placed on $V$, on $W$, and on the
cross edges \cite{cqsupper,cqslower}; fortunately here we shall see that the structure on $V$
suffices to determine the lower bounds.  Because every $N(v+w,v;4,C')$ is an $N(v+w,v;4,C'+1)$ for
$1 \leq C' \leq 3$, and $N(v+w,v;4,4)$ coincides with $N(v+w,4)$, $cost\ \mathscr{ON}(v+w,v;4,1)
\geq cost\ \mathscr{ON}(v+w,v;4,2) \geq cost\ \mathscr{ON}(v+w,v;4,3) \geq cost \
\mathscr{ON}(v+w,4)$. We use these obvious facts to establish lower and upper bounds without
further comment.

\section{Case $C'=1$}\label{sec:1}

\subsection{$\mathscr{ON}(n,v;4,1)$}

\begin{thm}\label{41} Let $n=v+w \geq 5$.
\begin{enumerate}
\item
$cost\ \mathscr{ON}(v+w,v;4,1) = cost \ \mathscr{ON}(v+w,4)$ when $v \leq w+1$.  \item $cost\
\mathscr{ON}(v+w,v;4,1) = \binom{v+w}{2}+\binom{v}{2} - \lfloor \frac{vw}{2} \rfloor$ when $v \geq
w+1$.
\end{enumerate}
\end{thm}
\begin{proof} To prove the lower
bound, we establish that $cost\ \mathscr{ON}(v+w,v;4,1) \geq \binom{v+w}{2}+\binom{v}{2} - \lfloor
\frac{vw}{2} \rfloor$.  It suffices to prove that the number of subgraphs employed in an
$N(v+w,v;4,1)$ other than triangles, kites, and 4-cycles is at least $ \lceil \binom{v}{2} -
\frac{1}{2}vw \rceil = \binom{v}{2} - \lfloor \frac{1}{2}vw \rfloor$.  By Lemma \ref{neutral}, this
is a lower bound on the number of positive edges in any such decomposition; because each positive
edge lies in a different subgraph of the decomposition, the lower bound follows.

Now we turn to the upper bounds. For the first statement, because an $\mathscr{ON}(v+w,v;4,1)$ is
also an $\mathscr{ON}(v+w,v-1;4,1)$, it suffices to consider $v\in\{w,w+1\}$. When $v=w$, write $v
= 4s+t$ with $t \in \{0,3,5,6\}$. Form on $V$ a complete multipartite graph with $s$ classes of
size  four and one class of size $t$. Replace  edge $e=\{x,y\}$ of this graph by the 4-cycle
$(x,y,a_x,a_y)$. On $\{x_1,\dots,x_{\ell},a_{x_1},\dots,a_{x_\ell}\}$ whenever
$\{x_1,\dots,x_\ell\}$ forms a class of the multipartite graph, place a decomposition  that is
optimal for drop cost and uses 4, 7, 12, and 17 wavelengths when $\ell$ is 3, 4, 5, or 6,
respectively (see Appendix \ref{ap:1}).

Now let $v=w+1$. Let $V=\{0,\dots,v-1\}$ and $W=\{a_0,\dots,a_{v-2}\}$.  Form triangles
$(i,i+1,a_i)$ for $0\leq i < v-1$.  Then form 4-cycles $(i,j+1,a_i,a_j)$ for $0 \leq i < j \leq
v-2$.

Finally, suppose that $v \geq w+2$. When $v$ is even, form a 1-factorization $F_0,\dots,F_{v-2}$ on
$V$.  For $0 \leq i < w$, let $\{ e_{ij} : 1 \leq j \leq \frac{v}{2}\}$ be the edges of $F_i$, and
form triangles $T_{ij} = \{a_i\} \cup e_{ij}$. Now for $0 \leq i < w$; $1 \leq j \leq \lfloor
\frac{w}{2} \rfloor$;  and furthermore $j \neq
 \frac{w}{2}$ if $i \geq  \frac{w}{2}$ and
$w$ is even, adjoin edge $\{a_i,a_{i+j\bmod w}\}$ to $T_{ij}$ to form a kite. All edges of
1-factors $\{F_i: w \leq i < v-1\}$  are taken as $K_2$s.

When $v$ is odd, form a near 1-factorization $F_0,\dots,F_{v-1}$ on $V$, in which $F_{v-1}$
contains the edges  $\{\{2h,2h+1\} : 0\leq h < \frac{v-1}{2}\}$, and near 1-factor $F_i$ misses
vertex $i$ for $0\leq i < v$.  Then form 4-cycles $(2h,2h+1,a_{2h+1},a_{2h})$ for $0 \leq h <
\lfloor \frac{w}{2} \rfloor$. For $0 \leq i < w$, let $\{ e_{ij} : 1 \leq j \leq \frac{v-1}{2}\}$
be the edges of $F_i$, and form  triangles $T_{ij} = \{a_i\} \cup e_{ij}$. Without loss of
generality we assume that $w-1 \in e_{01}$; when $w$ is odd, adjoin $\{w-1, a_{w-1}\}$ to $T_{01}$
to form a kite. Now for $0 \leq i < w$; $1 \leq j \leq \lfloor \frac{w}{2} \rfloor$; and
furthermore $j \neq
 \frac{w}{2}$ if $i \geq  \frac{w}{2}$ and $w$ is even and $j \neq 1$
if $i = 2h $ for $0 \leq h < \lfloor \frac{w}{2} \rfloor$, adjoin edge $\{a_i,a_{i+j\bmod w}\}$ to
$T_{ij}$ to form a kite. All edges of near 1-factors $\{F_i: w \leq i < v-1\}$ and the
$\frac{v-1}{2} - \lfloor \frac{w}{2} \rfloor$ remaining edges of $F_{v-1}$ are taken as $K_2$s.

When $v \geq w+1$, each subgraph contains exactly one edge on $V$ and so their number is ${v
\choose 2}$. This fact is later used to prove Theorem \ref{mon41}.
\end{proof}

\subsection{$\mathscr{MON}(n,v;4,1)$}

\begin{thm} \label{mon41small}
Let $v+w \geq 5$.  For $C'=1$ and $v \leq w$,
\[ wavecost \, \mathscr{MON}(v+w,v;4,1) = wavecost \, \mathscr{MON}(v+w,4). \]
\end{thm}
\begin{proof}
We need only treat the cases when $v \in \{w,w-1\}$; the case with $v=w$ is handled in the proof of
Theorem \ref{41}.  When $v=w-1$, the argument is identical to that proof, except that we choose
$v=4s+t$ with $t \in \{0,1,2,3\}$ and place decompositions on
$\{x_1,\dots,x_{\ell},a_{x_1},\dots,a_{x_\ell},a_v\}$ instead, with 1,3,6,9 wavelengths when $\ell
= 1,2,3,4$ respectively (see Appendix \ref{ap:2}).

\end{proof}

\begin{thm}\label{mon41}When $v > w$,
\[ wavecost \, \mathscr{MON}(v+w,v;4,1) = \binom{v}{2}. \]
\end{thm}
\begin{proof}
  Since  every edge on $V$ appears on a different wavelength, $\binom{v}{2}$
  is a lower bound.  As noted in the proof of Theorem
  \ref{41} the constructions given there meet this bound.
\end{proof}

The solutions used from Theorem \ref{41} are (essentially) the only ones to minimize the number of
graphs in an $\mathscr{ON}(v+w,v;4,1)$ with $v > w$.  However, perhaps surprisingly they are not
the only ones to minimize the number of wavelengths.  To see this, consider a
$\mathscr{ON}(v+w,v;4,1)$ with $v > w > 2$ from Theorem \ref{41}. Remove edges $\{a_0,a_1\}$,
$\{a_0,a_2\}$, and $\{a_1,a_2\}$ from their kites, and form a triangle from them.  This does not
change the drop cost, so the result is also an $\mathscr{ON}(v+w,v;4,1)$.  It has one more graph
than the original.  Despite this, it does not need an additional wavelength, since the triangle
$(a_0,a_1,a_2)$ can share a wavelength with an edge on $V$.  In this case, while minimizing the
number of connected graphs serves to minimize the number of wavelengths, it is not the only way to
do so.

\section{Case $C'=2$}\label{sec:2}

\subsection{$\mathscr{ON}(n,v;4,2)$}

\begin{thm}\label{42e} Let $v+w \geq 5$ and $v$ be even.
\begin{enumerate} \item When $v \leq 2w$, $cost\
\mathscr{ON}(v+w,v;4,2) = cost \ \mathscr{ON}(v+w,4)$.  \item When $v \geq 2w+2$, $cost\
\mathscr{ON}(v+w,v;4,2) = \binom{v+w}{2}+ \lceil \frac{1}{2}\binom{v}{2} \rceil - \frac{vw}{2} +
\delta$, where $\delta = 1$ if $w=4$ or if $w=2$ and $v \equiv 0 \pmod{4}$, and $\delta=0$
otherwise.
\end{enumerate}
 \end{thm}
\begin{proof}
By Lemma \ref{neutral}, $ \binom{v}{2} - vw$ is a lower bound on the number of positive edges in
any $N(v+w,v;4,2)$; every subgraph of the decomposition containing a positive edge contains at most
two positive edges. So the number of subgraphs employed in an  $N(v+w,v;4,2)$
 other than triangles, kites, and 4-cycles is at least
$\lceil \frac{1}{2} \left ( \binom{v}{2} - vw \right ) \rceil$. The lower bound follows for $w \neq
2,4$.

As in the proof of Lemma \ref{neutral},  denote by
 $\alpha$ (resp.  $\beta$) the number of subgraphs containing $2$
(resp $1$) neutral edges and so at least two cross edges. We have $2\alpha+\beta \leq
2\alpha+2\beta \leq vw $.
 Equality in the lower bound, when $v\equiv 0
\pmod{4}$, arises only when $\beta = 0$ and therefore to meet the bound an $\mathscr{ON}(w,4)$ must
be placed on $W$ implying that $\delta = 1$ if $ w = 2$ or $4$. When $v \equiv 2 \pmod{4}$, we can
have $2\alpha+\beta = vw -1$ and so  $\beta = 1$. We can use
 an edge on $W$ in  a graph with an
edge on $V$. But when $w=4$, the five edges that would remain on $W$ require drop cost 6, and so
$\delta = 1$.

Now we turn to the upper bounds. If $w \geq v-1$, apply Theorem \ref{41}. Suppose that $w \leq
v-2$. Let $V = \{0,\dots,2t-1\}$ and $W=\{a_0,\dots,a_{w-1}\}$.  Place an $ \mathscr{ON}(w,4)$ on
$W$. Form a 1-factorization on $V$ containing factors $ \{F_0,\dots,F_{w-1},G_0,\dots,G_{2t-2-w}\}$
in which the last two 1-factors are  $\{\{2h,2h+1\} : 0 \leq h < t\}$ and $\{\{2h+1,2h+2 \bmod 2t\}
: 0 \leq h < t\}$, whose union is a Hamilton cycle.  For $0 \leq i < w$, form triangles $T_{ij}$ by
adding $a_i$ to each edge $e_{ij} \in F_i$. For $0 \leq i < \mbox{min}(w,2t-1-w)$, observe that
$H_i = F_i \cup G_i$ is a 2-factor containing even cycles.  Hence there is a bijection $\sigma$
mapping edges of $F_i$ to edges of $G_i$ so that $e$ and $\sigma(e)$ share a vertex.  Adjoin edge
$\sigma(e_{ij})$ to the triangle $T_{ij}$ to form a kite.  In this way, all edges between $V$ and
$W$ appear in triangles or kites, and all edges on $V$ are employed when $v \leq 2w$. When $v \geq
2w+2$, the edges remaining on $V$ are those of the factors $G_{w},\dots,G_{v-2-w}$.

 When $v\neq
2w+2$, the union of these edges is connected because the union of the last two is connected, and
hence it can be partitioned into $P_3$s (and one $P_2$ when $v \equiv 2 \pmod{4}$)
\cite{CarSch,Yav}.  When $w=2$ and $v \equiv 2 \pmod{4}$, the drop cost can be reduced by 1 as
follows. Let $\{x,y\}$ be the $P_2$ in the decomposition, and let $\{x,z\} \in G_0$. Let $T$ be the
triangle obtained by removing $\{x,z\}$ from its kite.  Add $\{a_0,a_1\}$ to $T$ to form a kite.
Add the $P_3$ $[y,x,z]$. In this way two isolated $P_2$s are replaced by a $P_3$, lowering the drop
cost by 1.

When $v=2w+2$, we use a variant of this construction.  Let $R$ be a graph with vertex set $V$ that
is isomorphic to $\frac{v}{4}$ $K_4$s when $v \equiv 0 \pmod{4}$ and to $\frac{v-6}{4}$ $K_4$s and
one $K_{3,3}$ when $v \equiv 2 \pmod{4}$.  Let $F_1,\dots,F_{w-1},G_1,\dots,G_{w-1}$ be the
1-factors of a 1-factorization of the complement of $R$ (one always exists \cite{RosaWallis}).
Proceed as above to form kites using $a_i$ for $1 \leq i < w$ and the edges of $F_i$ and $G_i$.
For each $K_4$ of $R$ with vertices $\{p,q,r,s\}$, form kites $(a_0,q,p;r)$ and $(a_0,r,s;p)$.
Then add the $P_3$ $[r,q,s]$.  If $R$ contains a $K_{3,3}$ with bipartition
$\{\{p,q,r\},\{s,t,u\}\}$, add kites $(a_0,s,p;t)$, $(a_0,q,t;r)$, and $(a_0,r,u;p)$.  What remains
is the $P_4$ $[r,s;q,u]$, which can be partitioned into a $P_2$ and a $P_3$.
\end{proof}

In order to treat the odd case, we establish an easy preliminary result:

\begin{lemma}\label{cocktail} Let $w>3$ be a positive integer.
The graph on $w$ vertices containing all edges except for $\lfloor \frac{w}{2}\rfloor$ disjoint
edges (i.e., $K_w \setminus \lfloor \frac{w}{2} \rfloor K_2 $) can be partitioned into
\begin{enumerate}
 \item 4-cycles when $w$ is even;
\item kites and
4-cycles when $w \equiv 1 \pmod{4}$; and
 \item kites, 4-cycles, and
exactly two triangles when $w \equiv 3 \pmod{4}$.
\end{enumerate}
\end{lemma}
\begin{proof} Let $W =
\{a_0,\dots,a_{w-1}\}$. When $w$ is even, form 4-cycles $\{(a_{2i},a_{2j},a_{2i+1},a_{2j+1}): 0
\leq i < j < \frac{w}{2}\}$ leaving uncovered the $\frac{w}{2}$ edges $\{a_{2i}, a_{2i+1}\}$. (This
is also a consequence of a much more general result in \cite{FuR00}.)

 When $w$ is odd, the proof is by induction on $w$ by adding four new vertices.
 So we provide two base cases for the induction to cover all odd values of $w$.

For $w=5$, $K_5\setminus\{\{a_0,a_1\},\{a_2,a_3\}\}$ can be partitioned into the two kites
$(a_2,a_4,a_0;a_3)$ and $(a_3,a_4,a_1;a_2)$.

For $w=7$, $K_7\setminus\{\{a_0,a_1\},\{a_2,a_3\}, \{a_4,a_5\}\}$ can be partitioned into the kites
$(a_3,a_6,a_0;a_5)$, $(a_1,a_6,a_4; a_3)$ and $(a_5,a_6,a_2;a_1)$, and the $K_3$s $(a_0,a_2,a_4)$
and $(a_1,a_3,a_5)$.

By induction consider an optimal decomposition of $K_w-F$, with $F = \{\{a_{2h},a_{2h+1}\} : 0 \leq
h < \frac{w-1}{2}\}$.   Add four vertices $a_w,a_{w+1},a_{w+2},a_{w+3}$. Add the  $C_4$s
$(a_{2h},a_w,a_{2h+1},a_{w+1})$ and $ (a_{2h},a_{w+2},a_{2h+1},a_{w+3})$ where $0 \leq  h <
\frac{w-1}{2}$. Cover the edges of the $K_5$ on $\{a_{w-1},a_w,a_{w+1},a_{w+2},a_{w+3}\}$ minus the
edges $\{a_{w-1},a_w\}$ and $\{a_{w+1},a_{w+2}\}$, using two kites as shown for the case when
$w=5$.
\end{proof}

\begin{thm}\label{42o} Let $v+w \geq 5$ and $v$ be odd.
\begin{enumerate}
\item When $v \leq 2w-1$, $cost\
\mathscr{ON}(v+w,v;4,2) = cost \ \mathscr{ON}(v+w,4)$.
\item When $v
\geq 2w+1$, $cost\ \mathscr{ON}(v+w,v;4,2) = \binom{v+w}{2}+ \lceil \frac{1}{2} \left (
\binom{v}{2} - vw + \lceil \frac{w}{2} \rceil \right ) \rceil
 + \delta$, where $\delta = 1$ if $w=3$ and $v \equiv 3
\pmod{4}$, $0$ otherwise.
\end{enumerate}
\end{thm}
\begin{proof} To prove the lower bound, it
suffices to prove that the number of subgraphs employed in an $N(v+w,v;4,2)$ other than triangles,
kites, and 4-cycles is at least $\lceil \frac{1}{2} \left ( \binom{v}{2} - vw + \lceil \frac{w}{2}
\rceil \right ) \rceil$.  As in the proof of Theorem \ref{42e}, this follows from Lemma
\ref{neutral}. When $w =3$ and $v \equiv 3 \pmod{4}$, at least $\binom{v}{2} - 3v + 2$ edges are
positive, an even number.  To meet the bound, exactly one cross edge remains and exactly two edges
on $W$ remain. These necessitate a further graph that is not a triangle, kite, or 4-cycle.

Now we turn to the upper bounds.  By Theorem \ref{42e}, $cost\ \mathscr{ON}((v+1)+(w-1),v+1;4,2) =
cost \ \mathscr{ON}(v+w,4)$ when $v\leq 2w-3$.  So suppose that $v \geq 2w-1$.  Write $v=2t+1$.

When $w=t+1$, form a near 1-factorization on $V$ consisting of $2t+1$ near 1-factors,
$F_0,\dots,F_{t},$ $G_0,\dots,G_{t-1}$. Without loss of generality, $F_i$ misses vertex $i$ for $0
\leq i \leq t$, and $F_t$ contains the edges $\{\{k,t+k+1\} : 0 \leq k < t\}$. The union of any two
near 1-factors contains a nonnegative number of even cycles and a path with an even number of
edges. For $0 \leq i \leq t$, form triangles $T_{ij}$ by adding $a_i$ to each edge $e_{ij} \in
F_i$. As in the proof of Theorem \ref{42e}, for $0 \leq i < t$, use the edges of $G_i$ to convert
every triangle $T_{ij}$ into a kite.  Then add edge $\{i,a_i\}$ to triangle $T_{ti}$ constructed
from edge $\{i,t+1+i\}$. What remains is the single edge $\{t,a_t\}$ together with all edges on
$W$.

When $w \not\in\{2,4\}$, place an $\mathscr{ON}(w,4)$ on $W$ of cost $\binom{w}{2}$ so that $a_t$
appears in a triangle in the decomposition, and use the edge $\{t,a_t\}$ to convert this to a kite.
We use a decomposition having $1 \leq \delta \leq 4$ triangles, therefore  getting a solution with
at most 3 triangles. Such a decomposition exists by Theorem \ref{lavbermcoudhu} if $ w \nequiv
0,1\pmod{8}$. If $ w \equiv 0,1\pmod{8}$ we build  a solution using $4$ triangles as follows. If $w
\equiv 1 \pmod{8}$, form an $\mathscr{ON}(w-2,4)$ on vertices $\{0,\dots,w-3\}$ with $3$ triangles.
Add the triangle $(w-3,w-2,w-1)$ and the 4-cycles $\{(2h,w-2,2h+1,w-1) : 0 \leq h <
\frac{w-3}{2}\}$. For $w =8$ a solution with 4 triangles is given in Appendix \ref{ap:3}. In
general, for $w \equiv 0 \pmod{8}$, form an $\mathscr{ON}(w-8,4)$ on vertices $\{0,\dots,w-9\}$
with $4$ triangles. Add the 4-cycles $\{(2h,w-2j,2h+1,w-2j + 1) : 0 \leq h < \frac{w-8}{2}\}; 1\leq
j \leq 4$ and an $\mathscr{ON}(8,4)$ without
 triangles on the $8$ vertices $\{w-8,\dots,w-1\}$.

Two values for $w$ remain.  When $w=2$, an $\mathscr{ON}(5,3;4,1)$ is also an
$\mathscr{ON}(5,3;4,2)$. The case when $v=7$ and $w=4$ is given in Appendix \ref{ap:3}. The
solution given has only $1$ triangle.

Henceforth $w \leq t$.  For $t>2$, form a near 1-factorization
$\{F_0,\dots,F_{w-1},G_0,\dots,G_{2t-1-w}\}$ of $K_v\setminus C_t$, where $C_t$ is the $t$-cycle on
$(0,1,\dots,t-1)$; such a factorization exists \cite{Plantholt81}. Name the factors so that the
missing vertex in $F_i$ is $\lfloor i/2 \rfloor$ for $0 \leq i < w$ (this can be done, as every
vertex $i$ satisfying $0 \leq  i < t$ is the missing vertex in two of the near 1-factors).  Form
triangles using $F_0,\dots,F_{w-1}$ and convert to kites using $G_0,\dots,G_{w-1}$ as before.
There remain $2(t-w)$ near 1-factors $G_w, \dots, G_{2t-1-w}$.  For $ 0 \leq h < t-w$, $G_{w+2h}
\cup G_{w+2h+1}$ contains even cycles and an even path, and so partitions into $P_3$s.  Then the
edges remaining are (1) the edges of the $t$-cycle; (2) the edges $\{\{\lfloor i/2 \rfloor,a_i\}:0
\leq i < w\}$; and (3) all edges on $W$.  For $0 \leq i < \lfloor \frac{w}{2} \rfloor$, form
triangle $(i,a_{2i},a_{2i+1})$ and add edge $\{i,i+1\}$ to convert it to a kite.  Edges $\{\{i,i+1
\bmod t\} : \lfloor \frac{w}{2} \rfloor \leq i < t\}$ of the cycle remain from (1); edge $\{
\frac{w-1}{2}, a_{w-1} \}$ remains when $w$ is odd, and no edge remains when $w$ is even, from (2);
and all edges excepting a set of $\lfloor \frac{w}{2} \rfloor$ disjoint edges on $W$ remain.

When $w\neq 3$, we partition the remaining edges in (1) (which form a path of length $t-\lfloor
\frac{w}{2} \rfloor$),  into $P_3$s when $t-\lfloor \frac{w}{2} \rfloor$ is even, and into $P_3$s
and the $P_2$ $\{0,t-1\}$ when $t-\lfloor \frac{w}{2} \rfloor$ is odd. We adjoin edge $\{
\frac{w-1}{2}, a_{w-1} \}$ to the $P_3$ (from the $t$-cycle) containing the vertex $\frac{w-1}{2}$
to form a $P_4$. Finally, we  apply Lemma \ref{cocktail} to exhaust the remaining edges on $W$.

When $w=3$, the remaining edges are those of the  path $[0,t-1,t-2,\dots,2,1,a_2]$ and edges $\{
\{a_{2},a_0\}, \{a_{2},a_1\}\}$.  Include $\{ \{1,2\}, \{1,a_{2}\}, \{a_{2},a_0\}, \{a_{2},a_1\}\}$
in the decomposition, and partition the remainder into $P_3$s and, when $v \equiv 3 \pmod{4}$, one
$P_2$ $\{0,t-1\}$.

The case when $t=2$ is done in Example \ref{ex:2} (the construction is exactly that given above,
except that we start with a near 1-factorization of $K_5\setminus\{\{0,1\},\{0,3\}\}$).
\end{proof}

\subsection{$\mathscr{MON}(n,v;4,2)$}

\begin{thm} \label{mon42small}
For $C'=2$ and $v \leq 2w$,
\[ wavecost \, \mathscr{MON}(v+w,v;4,2) = wavecost \, \mathscr{MON}(v+w,4) . \]
\end{thm}
\begin{proof}
  It suffices to prove the statement for $v \in
  \{  2w-2,2w-1,2w\}$.  When $v=2w-1$, apply the construction given in the proof of Theorem
  \ref{42o},  where we noted that there are at most 3 triangles.
The  proof of Theorem \ref{42o} provides explicit solutions when $w \in  \{2,4\}$.

  Now suppose that $v=2w$.  In the proof of Theorem \ref{42e},
  $\frac{v}{2}=w$ triangles containing one edge on $V$ and two edges
  between a vertex of $V$ and $a_{w-1}$ remain.  Then convert $w-1$
  triangles to kites using edges on $W$ incident to $a_{w-1}$. That leaves one triangle. When
  the remaining edges on the $w-1$ vertices of $W$ support a
  $\mathscr{MON}(w-1,4)$ that contains at most two triangles, we are
  done.  It remains to treat the cases when $w-1 \equiv 2,7 \pmod{8}$
  or $w-1=4$.For the first case, let $x$ be one vertex of the triangle left
containing $a_{w-1}$, namely $(a_{w-1},x,y)$. Consider the pendant edge $\{x,t\} \in
  G_{w-2}$ used in a  kite containing $a_{w-2}$. Delete  $\{x,t\}$
 from this kite and adjoin $\{a_{w-3},a_{w-2}\}$ to
  the unique triangle so formed forming another kite.
 Finally adjoin $\{x,t\}$  to the triangle $(a_{w-1},x,y)$.
  Proceed as before, but partition all edges on
  $\{a_0,\dots,a_{w-2}\}$ except edge $\{a_{w-3},a_{w-2}\}$ into
  4-cycles and kites. The case when $w-1 = 4$ is similar, but we leave
  three of the triangles arising from $F_{w-1}$ and partition
  $K_5\setminus P_3$ into two kites.

 Now suppose that $v=2w- 2$. We do a construction similar to that above.
In the proof of Theorem \ref{42e}, there remain
  $3 \frac{v}{2}=3(w-1)$ triangles joining $a_{w-3}$ (resp.
  $a_{w-2},a_{w-1}$)
to $F_{w-3}$ (resp. $F_{w-2},F_{w-1}$). Then convert the $w-1$
  triangles containing $a_{w-1}$ to kites using edges
on $W$ incident to $a_{w-1}$,  $w-2$
  triangles containing $a_{w-2}$ to kites using the remaining edges
on $W$ incident to $a_{w-2}$, and $w-3$
  triangles containing $a_{w-3}$ to kites using edges
on $W$ incident to $a_{w-3}$.
 That leaves three triangles. So, if $w-3 \equiv 0,1 \pmod{8}$ we are done.
Otherwise, as above, choose in each of the three remaining triangles vertices $x_1,x_2,x_3$;
consider the edges $\{x_1,t_1\}$ (resp. $\{x_2,t_2\}$)
  appearing in the kites containing $a_{w-4}$ and $x_1$ (resp. $a_{w-4}$
and $x_2$), and the edge $\{x_3,t_3\}$ in the kite  containing $a_{w-5}$ and $x_3$. Delete these
edges and adjoin them to the three remaining triangles. Finally adjoin the edges
$\{a_{w-4},a_{w-5}\}$ and $\{a_{w-4},a_{w-6}\}$ to the two triangles
 obtained from the two kites  containing  $a_{w-4}$,
 and adjoin the edge $\{a_{w-5},a_{w-6}\}$ to the triangle
 obtained from the kite  containing  $a_{w-5}$.
 Proceed as before, but partition all edges on
  $\{a_0,\dots,a_{w-4}\}$ except the triangle $(a_{w-6},a_{w-5},a_{w-4})$
into 4-cycles and kites.

\end{proof}

\begin{thm} \label{mon42} \begin{enumerate} \item When $v > 2w$ is even,
\[ wavecost \, \mathscr{MON}(v+w,v;4,2) =
\left  \lceil  \left ( {2\binom{v}{2} + \binom{w}{2} } \right ) / 4 \right \rceil . \]
\item When $v > 2w$ is odd,
\[ wavecost \, \mathscr{MON}(v+w,v;4,2) = \left  \lceil  \left
( {2\binom{v}{2} + \frac{(w-1)(w+1)}{2} }  \right ) / 4 \right \rceil . \]
\end{enumerate}
\end{thm}
\begin{proof}
  First we treat the case when $v$ is even.  Then (by Theorem
  \ref{42e}) an $\mathscr{ON}(v+w,v;4,2)$ must employ $vw$ or $vw-1$
  neutral edges, using all $vw$ edges between $V$ and $W$.  Each such
  graph uses two edges on $V$ and none on $W$, except that a single
  graph may use one on $V$ and one on $W$.  Now the edges of $V$ must
  appear on $\lceil \frac{1}{2} \binom{v}{2} \rceil$ different
  wavelengths, and these wavelengths use at most one edge on $W$ (when
  $v \equiv 2 \pmod{4}$). Thus at least $\lceil \binom{w}{2}/4 \rceil$
  additional wavelengths are needed when $v \equiv 0 \pmod{4}$, for a
  total of $\lceil \binom{v}{2}/2 + \binom{w}{2}/4 \rceil$.  When $v
  \equiv 2 \pmod{4}$, at least $\lceil (\binom{w}{2}-1)/4 \rceil$
  additional wavelengths are needed; again the total is $\lceil
  \binom{v}{2}/2 + \binom{w}{2}/4 \rceil$.  Theorem \ref{42e} realizes
  this bound.

  When $v$ is odd, first suppose that $w$ is even. In order to realize
  the bound of Theorem \ref{42o} for drop cost,  by Lemma \ref{neutral},
$\frac{w}{2}$ neutral
  edges appear in subgraphs with one neutral edge and all other
  neutral edges appear in subgraphs with two.  In both cases, two
  edges between $V$ and $W$ are consumed by such a subgraph. When two
  neutral edges are used, no edge on $W$ can be used ; when one neutral edge
  is used, one edge on $W$ can also be used. It follows that the
  number of wavelengths is at least $\frac{1}{2} (\binom{v}{2} -
  \frac{w}{2}) + \frac{w}{2} + \frac{1}{4}(\binom{w}{2} -
  \frac{w}{2})$. This establishes the lower bound.  The case when $w$
  is odd is similar. The proof of Theorem \ref{42o} gives constructions with at most $3$ triangles and so establishes the upper bound except when
  $v \equiv 1 \pmod{4}$ and $w \equiv 3 \pmod{4}$, $w \neq 3$, where
  the construction employs one more graph than the number of
  wavelengths permitted. However, one graph included is the $P_2$
  $\{0,t-1\}$, and in the decomposition on $W$, there is a triangle.
  These can be placed on the same wavelength to realize the bound.
\end{proof}

When $v \equiv 1 \pmod{4}$ and $w \equiv 3 \pmod{4}$, $w \neq 3$, we place a disconnected graph,
$P_2 \cup K_3$, on one wavelength in order to meet the bound.  The construction of Theorem
\ref{42o} could be modified to avoid this by instead using a decomposition of $K_w \setminus (K_3
\cup \frac{w-3}{2} K_2)$ into 4-cycles and kites, and using the strategy used in the case for
$w=3$.  In this way, one could prove the slightly stronger result that the number of (connected)
subgraphs in the decomposition matches the lower bound on number of wavelengths needed.

In Theorem \ref{mon41}, the number of wavelengths and the drop cost are minimized simultaneously by
the constructions given; each constructed $\mathscr{ON}(v+w,v;4,1)$ has not only the minimum drop
cost but also the minimum number of wavelengths over all $N(v+w,v;4,1)$s. This is not the case in
Theorem \ref{mon42}.  For example, when $v > (1+\sqrt{2})w$, it is easy to construct an
$N(v+w,v;4,2)$ that employs only $\lceil \binom{v}{2} / 2 \rceil$ wavelengths, which is often much
less than are used in Theorem \ref{mon42}.  We emphasize therefore that a
$\mathscr{MON}(v+w,v;4,2)$ minimizes the number of wavelengths over all $\mathscr{ON}(v+w,v;4,2)$s,
{\sl not necessarily} over all $N(v+w,v;4,2)$s.

\section{Case $C'=3$}\label{sec:3}

\subsection{$\mathscr{ON}(n,v;4,3)$}

\begin{thm} \label{43} Let $v+w \geq 5$.
\begin{enumerate}
\item When $w \geq 1$, $cost\ \mathscr{ON}(v+w,v;4,3) = cost \ \mathscr{ON}(v+w,4)$.
\item $cost\ \mathscr{ON}(v+0,v;4,3) = cost \ \mathscr{ON}(v,3)$.
\end{enumerate}
\end{thm}
\begin{proof}
The second statement is trivial. Moreover $cost \ \mathscr{ON}(n,4) = cost \ \mathscr{ON}(n,3)$
when $n\equiv 1,3 \pmod{6}$, and hence the first statement holds when $v+w \equiv 1,3 \pmod{6}$.
To complete the proof  it suffices to treat the upper bound when $w=1$.

When $v+1 \equiv 5 \pmod{6}$, there  is a maximal partial triple system $(X,{\cal B})$ with
$|X|=v+1$ covering all edges except those in the 4-cycle $(r,x,y,z)$.  Set $W=\{r\}$, $V=X\setminus
W$, and add the 4-cycle to the decomposition to obtain an $\mathscr{ON}(v+1,v;4,3)$.

When $v \equiv 1,5 \pmod{6}$, set $\ell=v-1$ and when $v\equiv 3 \pmod{6}$ set $\ell=v-3$. Then
$\ell$ is even.  Form a maximal partial triple system $(V,{\cal B})$, $|V|=v$, covering all edges
except those in an $\ell$-cycle $(0,1,\dots,\ell-1)$ \cite{ColbournR86}.   Add a vertex $a_0$ and
form kites $(a_0,2i,2i+1;(2i+2)\bmod \ell)$ for $0 \leq i < \frac{\ell}{2}$.  For $i \in
\{\ell,\dots,v-1\}$, choose a triple $B_i \in {\cal B}$ so that $i \in B_i$ and $B_i = B_j$ only if
$i=j$.  Add $\{a_0,i\}$ to $B_i$ to form a kite. This yields an $\mathscr{ON}(v+1,v;4,3)$.
\end{proof}

\subsection{$\mathscr{MON}(n,v;4,3)$}

We focus first on lower bounds in Section \ref{sec:lower} and then we provide constructions
attaining these lower bounds in Section \ref{sec:upper}.

\subsubsection{Lower Bounds}
\label{sec:lower}

When $C'=3$, Theorem \ref{43} makes no attempt to minimize the number of wavelengths.  We focus on
this case here. Except when $n \in \{2,4\}$ or $v = n$, $cost\, \mathscr{ON}(n,v;4,3) =
\binom{n}{2}$, and every graph in an $\mathscr{ON}(n,v;4,3)$ is a triangle, kite, or 4-cycle.  Let
$\delta$, $\kappa$, and $\gamma$ denote the numbers of triangles, kites, and 4-cycles in the
grooming, respectively. Then $3\delta + 4\kappa + 4\gamma = \binom{n}{2}$, and the number of
wavelengths is $\delta+\kappa+\gamma$. Thus in order to minimize the number of wavelengths, we must
minimize the number $\delta$ of triangles. We focus on this equivalent problem henceforth.

In an $\mathscr{ON}(n,v;4,3)$, for $0 \leq i \leq 3$ and $0 \leq j \leq 4$, let $\delta_{ij}$,
$\kappa_{ij}$, and $\gamma_{ij}$ denote the number of triangles, kites, and 4-cycles, respectively,
each having $i$ edges on $V$ and $j$ edges between $V$ and $W$. The only counts that can be nonzero
are $\delta_{00}$, $\delta_{02}$, $\delta_{12}$, $\delta_{30}$; $\kappa_{00}$, $\kappa_{01}$,
$\kappa_{02}$, $\kappa_{03}$, $\kappa_{12}$, $\kappa_{13}$, $\kappa_{22}$, $\kappa_{31}$;
$\gamma_{00}$, $\gamma_{02}$, $\gamma_{04}$, $\gamma_{12}$, $\gamma_{22}$. We write $\sigma_{ij} =
\kappa_{ij} + \gamma_{ij}$ when we do not need to distinguish kites and 4-cycles. Our objective is
to minimize $\delta_{00} + \delta_{02} + \delta_{12}+\delta_{30}$ subject to certain constraints;
we adopt the strategy of \cite{cqslower} and treat this as a linear program.

Let $\varepsilon = 0$ when $v \equiv 1,3 \pmod{6}$, $\varepsilon = 2$ when $v \equiv 5 \pmod{6}$,
and $\varepsilon = \frac{v}{2}$ when $v \equiv 0 \pmod{2}$. We specify the linear program in Figure
\ref{linprog}.  The first row lists the primal variables.  The second lists coefficients of the
objective function to be minimized.  The remainder list the coefficients of linear inequalities,
with the final column providing the {\sl lower bound} on the linear combination specified.  The
first inequality states that the number of edges on $V$  used is at least the total number on $V$,
while the second specifies that the number of edges used between $V$ and $W$ is at most the total
number between $V$ and $W$. For the third, when $v \equiv 5 \pmod{6}$ at least  four edges on $V$
are not in triangles, and so at least two graphs containing edges of $V$ do not have a triangle on
$V$; when $v \equiv 0 \pmod{2}$ every graph can induce at most two odd degree vertices on  $V$, yet
all are odd in the decomposition.

\begin{figure}[htbp]
\begin{center}
\begin{tabular}{|cccc|ccccccccc|c|}
\hline
$\delta_{30}$ & $\delta_{12}$ & $\delta_{02}$ & $\delta_{00}$ & $\kappa_{31}$ & $\sigma_{22}$ & $\kappa_{13}$& $\sigma_{12}$ & $\gamma_{04}$ & $\kappa_{03}$ & $\sigma_{02}$ &$\kappa_{01}$ & $\sigma_{00}$ & \\
\hline
1 & 1 & 1 & 1 & 0 & 0 & 0 & 0 & 0 & 0 & 0 & 0 & 0 & \\
\hline
3 & 1 & 0 & 0 & 3 & 2 & 1 & 1 & 0 & 0 & 0 & 0 & 0 & $\binom{v}{2}$\\
0 & -2 & -2 & 0 & -1 & -2 & -3 & -2 & -4 & -3 & -2 & -1 & 0 & $-vw$\\
0 & 1 & 0 & 0 & 0 & 1 & 1 & 1 & 0 & 0 & 0 & 0 & 0 & $\varepsilon$\\
\hline
\end{tabular}
\end{center}
\caption{The linear program for $\mathscr{ON}(n,v;4,3)$.} \label{linprog}
\end{figure}

We do not solve this linear program.  Rather we derive  lower bounds by considering its dual.  Let
$y_1$, $y_2$, and $y_3$ be the dual variables. A dual feasible solution has $y_1 = \frac{1}{3}$,
$y_2=1$, and $y_3=\frac{4}{3}$, yielding a dual objective function value of $\frac{1}{6}v(v-1) - vw
+ \frac{4}{3} \varepsilon$.  Recall that every dual feasible solution gives a lower bound on all
primal feasible solutions

On the other hand, $3 \delta \equiv \binom{n}{2} \pmod{4}$ and so $\delta \equiv 9\delta \equiv
3\binom{n}{2} \pmod{4}$. The value of $3\binom{n}{2} \pmod{4}$ is in fact the value of $t$ given in
Theorem \ref{lavbermcoudhu}. Therefore if $x$ is a lower bound on $\delta$ in an
$\mathscr{ON}(n,v;4,3)$, so is $\langle x \rangle_n$, where $\langle x \rangle_n$ denotes the
smallest nonnegative integer $\overline{x}$ such that $\overline{x} \geq x$ and $\overline{x}
\equiv 3\binom{n}{2} \pmod{4}$.

The discussion above proves the general lower bound on the number of triangles:

\begin{thm}\label{lower43}
Let $v+w \geq 5$, and let

\[ L(v,w) =  \left \{ \begin{array}{ccl}
\frac{1}{6}v(v-1) - vw & \mbox{if} & v \equiv 1,3 \pmod{6}\\
\frac{1}{6}v(v-1) - vw + \frac{8}{3} & \mbox{if} & v \equiv 5 \pmod{6}\\
\frac{1}{6}v(v+3) - vw  & \mbox{if} & v \equiv 0 \pmod{2}\\
\end{array} \right . \]

\noindent Then the  number of triangles in an $\mathscr{ON}(v+w,v;4,3)$ is at least
\[ \delta_{\min}(v,w) = \langle  L(v,w)  \rangle_{v+w} \]
\end{thm}

\begin{rmk}
\label{rmk:0} In particular, if $v$ is odd and  $w\geq \lceil \frac{v-1}{6}\rceil$ or if $v$ is
even and  $w\geq \lceil \frac{v-4}{6}\rceil$, then $L(v,w)\leq 0$ and the minimum number of
triangles is $\delta_{\min}(v,w) = \langle  0 \rangle_{v+w} \leq 3$.
\end{rmk}

\subsubsection{Upper Bounds}
\label{sec:upper}

We first state two simple lemmas to be used intensively in the proof of Theorem \ref{thm:MON43}.
The following result shows that in fact we do not need to check \emph{exactly} that the number of
triangles of an optimal construction meets the bound of Theorem \ref{lower43}.

\begin{lemma}
\label{lemma:interval} Any $\mathscr{ON}(v+w,v;4,3)$ is a $\mathscr{MON}(v+w,v;4,3)$ if the number
of triangles that it contains is at most $\max(3,\lceil L(v,w)\rceil+3)$.

\end{lemma}
\begin{proof}
In the closed interval $[\lceil L(v,w)\rceil,\lceil L(v,w)\rceil+3]$ there is exactly one integer
congruent to $3\binom{n}{2} \pmod{4}$, and so necessarily exactly one integer equal to
$\delta_{\min}(v,w)$.
\end{proof}

Combining Remark \ref{rmk:0} and Lemma \ref{lemma:interval} we deduce that when $v$ is odd and
$w\geq \lceil \frac{v-1}{6}\rceil$ or if $v$ is even and $w\geq \lceil \frac{v-4}{6}\rceil$, to
prove the optimality of a construction it is enough to check that there are at most three
triangles.

As a prelude to the constructions, let $(V,{\cal B})$ be a partial triple system, $V = \{0,\dots,
v-1\}$, and ${\cal B} = \{B_1,\dots,B_b\}$. Let $r_i$ be the number of blocks of ${\cal B}$ that
contain $i \in V$. A {\sl headset} is a multiset $S = \{s_1,\dots,s_b\}$ so that $s_k \in B_k$ for
$1 \leq k \leq b$, and for $0 \leq i \leq v-1$ the number of occurrences of $i$ in $S$ is $\lfloor
\frac{r_i}{3} \rfloor$ or $\lceil \frac{r_i}{3} \rceil$.

\begin{lemma}\label{equit}
Every partial triple system has a headset.
\end{lemma}
\begin{proof}
Form a bipartite graph $\Gamma$ with vertex set $V\cup{\cal B}$, and an edge $\{v,B\}$ for $v \in
V$ and $B \in {\cal B}$ if and only if $v \in B$. The graph $\Gamma$ admits an equitable
3-edge-colouring \cite{deW}; that is, the edges can be coloured green, white, and red so that every
vertex of degree $d$ is incident with either $\lfloor d/3\rfloor$ or $\lceil d/3 \rceil$ edges of
each colour. Then for $1 \leq k \leq b$,  $B_k$ is incident to exactly three edges, and hence to
exactly one edge $\{i_k,B_k\}$ that is green; set $s_k = i_k$. Then $(s_1,\dots,s_b)$ forms the
headset.
\end{proof}

\begin{thm}
\label{thm:MON43} Let $v+w \geq 5$. When $w \geq 1$,
\[ wavecost \, \mathscr{MON}(v+w,v;4,3) =
\left \lceil \left ( \binom{v+w}{2} + \delta_{\min}(v,w) \right ) / 4 \right  \rceil . \]
\end{thm}
\begin{proof}
The lower bound follows from Theorem \ref{lower43}, so we focus on the upper bound.

When $w \geq 1$, an $\mathscr{ON}(v+w,v;4,3)$ of cost $\binom{v+w}{2}$ is an
$\mathscr{ON}(v+w,v-1;4,3)$. Let us show that it suffices to prove the statement for $w \leq
\frac{v+9}{6}$ when $v$ is odd, and for $w \leq \frac{v+4}{6}$ when $v$ is even. Equivalently, we
show that if it is true for these values of $w$, then it follows for any $w$. Note that
$\delta_{\min}(v,w) \leq 3$ if $\delta_{\min}(v+1,w-1)\leq 3$.

Indeed, let $v$ be even. If $w=\lfloor \frac{v+4}{6}\rfloor+1$, the result follows from the case
for $v+1$ (odd) and $w-1=\lfloor \frac{v+4}{6}\rfloor \leq \frac{v+1+9}{6}$, in which case
$\delta_{\min}(v+1,w-1)=\langle 0\rangle_{v+w}$. If $w=\lfloor \frac{v+4}{6}\rfloor+2$ it follows
from the case for $v+1$ (odd) and $w-1=\lfloor \frac{v+4}{6}\rfloor + 1\leq \frac{v+1+9}{6}$, and
$\delta_{\min}(v+1,w-1)=\langle 0\rangle_{v+w}$. If $w \geq \lfloor \frac{v+4}{6}\rfloor+3$ it
follows from the case for $v+2$ (even) and $w-2$.

Let $v$ be odd. If $w=\lfloor \frac{v+9}{6}\rfloor+1$ it follows from the case for $v+1$ (even) and
$w-1$, which has been already proved (in this case also $\delta_{\min}(v+1,w-1)=\langle
0\rangle_{v+w}$). If $w \geq \lfloor \frac{v+9}{6}\rfloor+2$ it follows from the case for $v+2$
(odd) and $w-2$.

In each case, we use the same general prescription.  Given a partial triple system $(V,{\cal B})$,
a headset  $S = \{s_1,\dots,s_b\}$ is formed using Lemma \ref{equit}.  Add vertices
$W=\{a_0,\dots,a_{w-1}\}$, a set disjoint from $V$ of size $w \geq 1$. For each $i$ let $D_i$ be a
subset of $\{0,\ldots,w-1\}$, which is specified for each subcase, and that satisfies the following
property: $|D_i|$ is at most the number of occurrences of $i$ in the headset $S$. Among the blocks
$B_k$ such that $s_k=i$, we choose $|D_i|$ of them, namely the subset $\{B_k^j:j \in D_i\}$, and
form $|D_i|$ kites by adding for each $j \in D_i$ the edge $\{a_j,i\}$ to the block $B_k^j$.

The idea behind the construction is that if we can choose $|D_i|=w$, we use all the edges between
$V$ and $W$ leaving a minimum number of triangles in the partition of $V$ (see Case \textbf{O1a}).
Unfortunately it is not always possible to choose $|D_i|=w$, in particular when $w$ is greater than
the number of occurrences of $i$ in the headset.
So we distinguish different cases:\\

\noindent {\bf Case O1a. $v=6t+1$ or $6t+3$ and $w \leq \frac{v-1}{6}$.} Let $(V,{\cal B})$ be a
Steiner triple system. For $0 \leq i < v$, let $D_i=\{0,\ldots,w-1\}$. Apply the general
prescription. If $v=6t+1$, $i$ appears $t$ times in $S$ and $w \leq \frac{v-1}{6}=t$. If $v=6t+3$,
$i$ appears $t$ or $t+1$ times in $S$ and $w \leq t$. In both cases $|D_i|$ is at most the number
of occurrences of $i$ in $S$, so the construction applies and all the edges between $V$ and $W$ are
used in the kites. All the edges on $V$ are used and $\frac{v(v-1)}{6}-vw$ triangles remain.
Finally, it remains to partition the edges of $W$. When $w\not\in\{2,4\}$, form a
$\mathscr{MON}(w,4)$ on $W$, and doing so we have at most $\delta_{\min}$ triangles. If $w=2$ or
$w=4$ remove edges $\{a_0,0\}$ and $\{a_1,0\}$ from their kites and partition $K_W$ together
with these edges into a triangle ($w=2$) or two kites ($w=4$).\\

\noindent {\bf Case O1b. $v=6t+5$  and $w \leq \frac{v-1}{6}$.} Form a partial triple system
$(V,{\cal B})$  covering all edges except those in the $C_4$ $(0,1,2,3)$. For $0 \leq i \leq 3$,
let $D_i=\{0,\ldots,w-2\}$ and for $4 \leq i <v$ $D_i=\{0,\ldots,w-1\}$. Apply the general
prescription. Add the kites $(a_{w-1},1,2;3)$ and $(a_{w-1},3,0;1)$. Here again $i$ appears at
least $t$ times in $S$ and $w \leq t$. So $D_i$ is at most the number of occurrences of $i$ in $S$.
Again we have used all the edges on $V$ and all the edges between $V$ and $W$. It remains to
partition the edges of $W$, and this can be done as in
the Case {\bf O1a}.\\

\noindent {\bf Case O2. $v = 6t +3$ and  $w = t+1$, $v
> 3$.}  Form a partial triple system covering all edges except those on
the $v$-cycle $\{\{i,(i+1) \bmod v\} : 0 \leq i < v\}$ \cite{ColbournR86}. Set $D_i =
\{1,\dots,w-1\}$ for all $i$.  Apply the general prescription. Adjoin edges from $a_{0}$ to a
partition of the cycle, minus edge $\{0,v-1\}$, into $P_3$s. The only edge between $V$ and $W$ that
remains is $\{a_0,v-1\}$. When an $\mathscr{ON}(w,4)$ exists having 1, 2, 3, or 4 triangles, this
edge is used to convert a triangle to a kite.  This handles all cases except when $w \in \{2,4\}$.
In these cases, remove the pendant edge $\{a_1,v-1\}$ from its kite. When $w=2$, $\{a_0,a_1,v-1\}$
forms a triangle. When $w=4$,  partition the edges on $W$ together with
$\{a_0,v-1\}$ and $\{a_1,v-1\}$ into two kites.\\

\noindent {\bf Case O3. $v =6t+1$ and $w = t+1$.}

When $t=1$, a $\mathscr{MON}(7+2,7;4,3)$ has $\cB=\{(0,a_1,a_0;6)$, $(2,0,6;a_1)$, $ (3,0,4;a_1)$,
$(1,0,5;a_1)$, $(3,6,5;a_0)$, $ (4,6,1;a_1)$, $(3,2,1;a_0)$, $(5,2,4;a_0)$, $ (a_0,2,a_1,3)\}$.

A solution with $t=2$ is given in Appendix \ref{ap:4}.

When $t \geq 3$, form a 3-GDD of type $6^t$ with  groups $\{\{6p+q : 0 \leq q < 6\} : 0 \leq p <
t\}$.  Let $D_{6p+q} = \{0,\dots,w-2\} \setminus \{p\}$ for $0 \leq p < t$ and $0 \leq q < 6$.
Apply the general prescription. For $0 \leq p < t$,  on $\{6p+q : 0 \leq q < 6\} \cup \{v-1\} \cup
\{a_{w-1},a_p\}$ place a $\mathscr{MON}(7+2,7;4,3)$ obtained from the solution  $\cB$ for $t=1$, by
replacing $q$ by $6p+q : 0 \leq q < 6$, 6 by $v-1$ $a_0$ by $a_{w-1}$ and $a_1$ by $a_p$;  then
omit the kite  $(a_p,6p,a_{w-1};v-1)$.  All edges on $W$ remain; the edges $\{a_{w-1},6p\}$ and
$\{a_p,6p\}$ remain for $0 \leq p < t$, and the edge $\{a_{w-1},v-1\}$ remains.

Add the kites $(a_{w-2},6(w-2),a_{w-1};v-1)$ and for $0 \leq j < w-2 =t-1$
$(6j,a_{w-1},a_j;a_{w-2})$. If $w-2\not\in\{2,4\}$, that is $t\not\in\{3,5\}$, place a
$\mathscr{MON}(w-2,4)$ on $W-a_{w-2}-a_{w-1}$. Note that, as $3\binom{w-2}{2} \equiv
3\binom{v+w}{2} \pmod{4}$, we have the right number of triangles (at most $3$). If $w-2 \in\{2,4\}$
remove edges $\{a_0,w-2\}$ and $\{a_1,w-2\}$ from their kites, and partition
$K_w$ together with these edges.\\

\noindent {\bf Case O4. $v = 6t + 5$ and $w = t+1$.}

For $t=0$, a $\mathscr{MON}(5+1,5;4,3)$ has kites $(3,a_0,0;1)$, $(1,a_0,2;3)$, $(1,3,4;a_0)$, and
triangle $(0,2,4)$.

For $t=1$, let $V=\{0,\ldots,10\}$ and $W=\{a_0,a_1\}$. A $\mathscr{MON}(11+2,11;4,3)$ is formed by
using an $\mathscr{MON}(5+1,5;4,3)$ on $\{0,1,2,3,4\} \cup \{a_0\}$, and a partition of the
remaining edges, denoted by ${\cal Q}$, into 15 kites and a triangle. So we have two triangles,
attaining $\delta_{min}(11,2)$ as $13\equiv 5 {\pmod 8}$. The partition of ${\cal Q}$ is as
follows: the triangle $(a_0,a_1,10)$ and the kites $(0,6,5;a_0)$, $(1,8,6;a_0)$, $(2,9,7;a_0)$,
$(3,10,8;a_0)$, $(4,6,9;a_0)$, $(8,9,0;a_1)$, $(5,7,1;a_1)$, $(5,8,2;a_1)$, $(6,7,3;a_1)$,
$(5,10,4;a_1)$, $(3,9,5;a_1)$, $(2,10,6;a_1)$, $(0,10,7;a_1)$, $(4,7,8;a_1)$, and $(1,10,9;a_1)$.

For $t=2$, a $\mathscr{MON}(17+3,17;4,3)$ is given in Appendix \ref{ap:4}.

For $t \geq 3$, form a 3-GDD of type $6^t$ with  groups $\{\{6p+q : 0 \leq q < 6\} : 0 \leq p <
t\}$.  Let $D_{6p+q} = \{0,\dots,w-2\} \setminus \{p\}$ for $0 \leq p < t$ and $0 \leq q < 6$.
Apply the general prescription. There remain uncovered for each $p$ the edges of the set ${\cal
Q}_p$ obtained from the complete graph on the set of vertices $\{6p+q: 0 \leq q < 6\} \cup
\{v-5,v-4,v-3,v-2,v-1\} \cup \{a_{w-1},a_p\}$ minus the complete graph on $\{v-5,v-4,v-3,v-2,v-1\}
\cup \{a_{w-1}\}$.

To deal with the edges of ${\cal Q}_p$, we start from a partition of ${\cal Q}$, where we replace
pendant edges in kites as follows: Replace $\{a_1,4\}$ by $\{a_1,10\}$, $\{a_0,8\}$ by
$\{a_0,10\}$, and  $\{a_1,2\}$ by $\{a_0,8\}$. We delete the triangle $(a_0,a_1,10)$, resulting in
a new partition of ${\cal Q}$ into 15 kites and the 3 edges $\{a_0,a_1\}$, $\{a_1,2\}$, and
$\{a_1,4\}$. Then we obtain a partition of ${\cal Q}_p$ by replacing $\{0,1,2,3,4\}$ by
$\{v-5,v-4,v-3,v-2,v-1\}$, $q+5$ by $6p+q$ for $0 \leq q < 6$, $a_0$ by $a_{w-1}$, and $a_1$ by
$a_p$. At the end we get a partition of ${\cal Q}_p$ into 15 kites plus the 3 edges
$\{a_{w-1},a_p\}$, $\{a_p,v-3\}$, and $\{a_p,v-1\}$.

Now the $3t$ edges $\{\{a_{w-1},a_p\}, \{a_p,v-3\}, \{a_p,v-1\}: 0 \leq p < t\}$ plus the uncovered
edges of $K_W$ form a $K_{t+3}$ missing a triangle on $\{a_{w-1},v-3,v-1\}$.  If $t+3 \equiv
2,3,4,5,6,7 \pmod{8}$, use Theorem \ref{lavbermcoudhu} to form a $\mathscr{ON}(t+3,4)$ having a
triangle $(v-3,v-1,a_{w-1})$ and 0, 1, or 2 other triangles; remove the triangle
$(v-3,v-1,a_{w-1})$ to complete the solution with 1, 2, or 3 triangles (the triangle
$(v-5,v-3,v-1)$ is still present).  A variant is needed when $t+3 \equiv 0,1 \pmod{8}$.  In these
cases, form a $\mathscr{ON}(t+3,4)$ (having no triangles) in which $(v-3,a_{w-1},v-1;a_1)$ is a
kite.  Remove all edges of this  kite, and use edge $\{a_1,v-1\}$ to convert triangle
$(v-5,v-3,v-1)$ to a kite.

Finally, place a $\mathscr{MON}(5+1,5;4,3)$ on $\{v-5,v-4,v-3,v-2,v-1\} \cup \{a_0\}$. Altogether
we have a partition of all the edges using at most 3
triangles.\\

\noindent {\bf Case O5. $v = 6t+5$ and  $w = t+2$.}

When $t=0$,  partition all edges on $\{0,1,2,3,4\} \cup \{a_0,a_1\}$ except $\{a_0,a_1\}$ into
kites $(3,1,a_0;0)$, $ (3,2,a_1;0)$, $(a_1,1,4;2)$, $(0,1,2;a_0)$, and $(3,0,4;a_0)$.  Then a
$\mathscr{MON}(5+2,5;4,3)$ is obtained by removing pendant edges $\{a_0,0\}$ and $\{a_1,0\}$ and
adding triangle $(a_0,a_1,0)$.

When $t=1$, a $\mathscr{MON}(11+3,11;4,3)$ on $\{0,\dots,10\} \cup \{a_0,a_1,a_2\}$ is obtained by
taking the above partition on $\{0,1,2,3,4\} \cup \{a_0,a_1\}$, the triangle $(a_0,a_1,a_2)$, and a
partition of the remaining edges (which form a graph called $\cal Q$) into 11 kites and 6 4-cycles
as follows: kites $(2,9,7;a_0)$, $(4,5,10;a_0)$, $(2,10,6;a_1)$, $(4,6,9;a_2)$, $(7,10,0;a_2)$,
$(6,8,1;a_2)$, $(5,8,2;a_2)$, $(5,9,3;a_2)$, $(7,8,4;a_2)$, $(6,7,5;a_2)$, and $(9,10,8;a_1)$; and
4-cycles $(0,6,a_0,5)$, $(0,8,a_0,9)$, $(1,5,a_1,7)$, $(1,9,a_1,10)$, $(3,6,a_2,7)$, and
$(3,8,a_2,10)$.

A solution with $t=2$ is given in Appendix \ref{ap:4}.

When $t \geq 3$, form a 3-GDD of type $6^t$ with groups  $\{\{6p+q : 0 \leq q < 6\} : 0 \leq p <
t\}$.  Let $D_{6p+q} = \{0,\dots,w-3\} \setminus \{p\}$ for $0 \leq p < t$ and $0 \leq q < 6$.
Apply the general prescription. Add a partition of the complete graph on $\{v-5,v-4,v-3,v-2,v-1\}
\cup \{a_{w-2},a_{w-1}\}$ as in the case when $t=0$. It remains to partition, for each $p$, $0 \leq
p < t$, the graph ${\cal Q}_p$ is obtained from the complete graph on $\{6p+q : 0 \leq q < 6\} \cup
\{v-5,v-4,v-3,v-2,v-1\} \cup \{a_{w-2},a_{w-1},a_p\}$ minus the complete graph on
$\{v-5,v-4,v-3,v-2,v-1\} \cup \{a_{w-2},a_{w-1}\}$. This partition is obtained from that of $\cal
Q$ by replacing $\{0,1,2,3,4\}$ by $\{v-5,v-4,v-3,v-2,v-1\}$, $a_0$ by $a_{w-2}$, $a_1$ by
$a_{w-1}$, and $a_2$ by $a_p$. What remains is precisely the edges on $W$, so place a
$\mathscr{MON}(w,4)$ on $W$ to complete the
construction.\\

\noindent {\bf Case O6. $v = 6t+3$ and $w = t+2$.}

When $t=0$, a $\mathscr{MON}(3+2,3;4,3)$ has triangles $(a_0,0,1)$ and $\{a_1,1,2\}$ and 4-cycle
$(0,2,a_0,a_1)$.

When $t=1$, on $\{0,\dots,8\} \cup \{a_0,a_1,a_2\}$, place kites $(2,6,4;a_0)$, $(0,8,4;a_1)$,
$(0,5,7;a_1)$, $(3,6,0;a_2)$, $(1,7,4;a_2)$, $(5,8,2;a_2)$, $(1,6,5;a_2)$, $(2,7,3;a_2)$,
$(3,8,1;a_2)$, $(3,5,a_0;a_2)$, $(7,a_0,6;a_2)$, $(6,8,a_1;a_2)$, $(7,a_2,8;a_0)$, and 4-cycle
$(3,4,5,a_1)$. Adding the blocks of a $\mathscr{MON}(3+2,3;4,3)$ forms a
$\mathscr{MON}(9+3,9;4,3)$.

A solution with $t=2$ is given in Appendix \ref{ap:4}.

When $t \geq 3$, form a 3-GDD of type $6^t$ with  groups $\{\{6p+j : 0 \leq q < 6\} : 0 \leq p <
t\}$.  Let $D_{6p+q} = \{0,\dots,w-3\} \setminus \{p\}$ for $0 \leq p < t$ and $0 \leq q < 6$.
Apply the general prescription.  For $0 \leq p < t$,  on $\{6p+q : 0 \leq q < 6\} \cup
\{v-3,v-2,v-1\} \cup \{a_{w-2},a_{w-1},a_p\}$ place a $\mathscr{MON}(9+3,9;4,3)$, omitting a
$\mathscr{MON}(3+2,2;4,3)$ on $\{a_{w-2},a_{w-1},v-3,v-2,v-1\}$.  Place a
$\mathscr{MON}(3+2,2;4,3)$ on $\{a_{w-2},a_{w-1},v-3,v-2,v-1\}$.  Remove edges $\{a_0,a_{w-2}\}$
and $\{a_1,a_{w-1}\}$ from their kites, and convert the two triangles in the
$\mathscr{MON}(3+2,2;4,3)$ to kites using these. What remains is all edges on
$\{a_0,\dots,a_{w-3}\}$ and everything is in kites or 4-cycles excepting one triangle involving
$a_0$ and one involving $a_1$. If $w-2 \equiv 0,1,3,6 \pmod{8}$, place a $\mathscr{MON}(w-2,4)$ on
$\{a_0,\dots,a_{w-3}\}$.  Otherwise partition all edges on $\{a_0,\dots,a_{w-3}\}$ except
$\{a_0,a_2\}$ and $\{a_1,a_2\}$ into kites, 4-cycles, and at most one triangle, and use the last
two edges to form kites with the excess triangles involving $a_0$ and $a_1$.  The partition needed
is easily produced for $w-2 \in \{4,5,7,9\}$ and hence by induction
for all the required orders.\\

\noindent {\bf Case E1. $v \equiv 0 \pmod{2}$ and $w \leq \frac{v+2}{6}$.} Write $v=6t+s$ for $s
\in \{0,2,4\}$.  Let $L=(V,E)$ be a graph with edges
\[ \{\{3i,3i+1\},\{3i,3i+2\},\{3i+1,3i+2\} : 0 \leq i < t\} \cup
\{\{i,3t+i\}: 0 \leq i < 3t\}, \] together with $\{6t,6t+1\}$ when $s=2$ and with
$\{\{6t,6t+1\},\{6t,6t+2\},\{6t,6t+3\}\}$ when $s=4$. Let $(V,{\cal B})$ be a partial triple system
covering all edges except those in  $L$ (this is easily produced).   Let $D_i=\{0,\dots,w-2\}$ for
$0\leq i < v$. Apply the general prescription. For $0 \leq i < t$ and $j \in \{0,1,2\}$, form the
4-cycle $(a_{w-1},3i+((j+1) \bmod 3),3i+j,3t+3i+j)$. When $s=4$, form 4-cycle
$(a_{w-1},6t+2,6t,6t+3)$.  When $s \in \{2,4\}$, form a triangle $(a_{w-1},6t,6t+1)$. All edges on
$V$ are used and all edges on $W$ remain. All edges between $V$ and $W$ are used. Except when $w
\in\{2,4\}$,  or $w \equiv 2,7 \pmod{8}$ and $v \equiv 2,4 \pmod{6}$ form a $\mathscr{MON}(w,4)$ on
$W$ to complete the proof. When $w\equiv 2,7 \pmod{8}$ and $v \equiv 2,4 \pmod{6}$, convert
$\{a_{w-1},6t,6t+1\}$ to a kite using an edge of the $K_w$, and partition the $K_w\setminus K_2$
into kites and 4-cycles.  When $w \in \{2,4\}$, remove edges $\{a_0,0\}$ and $\{a_1,0\}$ from their
kites, and partition $K_w$ together with these
edges.\\

\noindent {\bf Case E2. $v \equiv 2 \pmod{6}$ and $w = \frac{v+4}{6}$.} Choose $m$ as large as
possible so that $m \leq \frac{v}{2}$, $m \leq \binom{w}{2}$, and $\binom{w}{2} - m \equiv 0
\pmod{4}$. Partition the $\binom{w}{2}$ edges on $W$ into sets $E_c$ and $E_o$ with $|E_c| = m$, so
that the edges on $E_o$ can be partitioned into kites and 4-cycles; this is easily done.  Place
these kites and 4-cycles on $W$.  Then let $\{e_i : 0 \leq i < m\}$ be the edges in $E_c$; let
$a_{f_i} \in e_i$ when $0 \leq i < m$; $f_i=0$ when $m \leq i < \frac{v-2}{2}$; and $f_{(v-2)/2} =
1$ if $m < \frac{v}{2}$. Next form a 3-GDD of type $2^{v/2}$ on $V$ so that $\{\{2i,2i+1\} : 0 \leq
i < \frac{v}{2}\}$ forms the groups, and ${\cal B}$ forms the blocks.    For $0 \leq i <
\frac{v}{2}$, let $D_{2i} = D_{2i+1} = \{0,\dots,w-1\} \setminus \{f_i\}$.    Apply the general
prescription.  Now for $0 \leq i < \frac{v}{2}$, form the triangle $(a_{f_i},2i,2i+1)$ and for $0
\leq i < m$ add edge $e_i$ to form a kite. At most three triangles remain {\em except when} $v \in
\{14,20\}$, where four triangles remain.  To treat these cases, we reduce the number of triangles;
without loss of generality,  the 3-GDD contains a triple $\{v-8,v-6,v-4\}$ in a kite with edge
$\{a_1,v-8\}$.  Remove this kite, and form kites $(a_0,v-7,v-8;v-6)$, $(a_0,v-5,v-6;v-4)$,
$(a_0,v-3,v-4;v-8)$,  and $(v-2,v-1,a_1;v-8)$.
\end{proof}

\begin{cor}\label{caso3wgrande} Let  $v\geq 4$ and $\mu_3(v)$ be
defined by:
\begin{center}
\begin{tabular}{|c|c|c|c|c|c|c|c|c|c|c|}
\hline $v$ & $6$ & $6t,t\geq 2$ & $1+6t$ & $2+6t$ & $9$ & $3+6t,t\geq 2$
& $4$ & $10$ & $4+6t,t\geq 2$ & $5+6t$ \\
\hline
$\mu_3(v)$  & $1$ & $1+t$ & $t$ & $1+t$ & $1$ & $1+t$ & $1$ & $2$ & $2+t$ & $1+t$ \\
\hline
\end{tabular}\vspace{4 mm}
\end{center}
Then  $wavecost \, \mathscr{MON}(v+w,v;4,3) =\left\lceil\frac{(v+w)(v+w-1)}{8}\right\rceil$ if and
only if $w\geq\mu_3(v)$.
\end{cor}

\section{Conclusions}

The determination of $cost \, \mathscr{ON}(n,v;\Ca,\Cb)$ appears to be easier when $C'=4$ than the
case for $C'=3$ settled in \cite{cqsupper,cqslower}.  Nevertheless the very flexibility in choosing
kites, 4-cycles, or triangles also results in a wide range of numbers of wavelengths among
decompositions with optimal drop cost.  This leads naturally to the question of minimizing the drop
cost and the number of wavelengths simultaneously.  In many cases, the minima for both can be
realized by a single decomposition.  However, it may happen that the two minimization criteria
compete.  Therefore we have determined the minimum number of wavelengths among all decompositions
of lowest drop cost for the specified values of $n$, $v$, and $C'$.

\section*{Acknowledgments}
Research of the authors is supported by MIUR-Italy (LG,GQ), European project IST FET AEOLUS, PACA
region of France, Ministerio de Educaci\'on y Ciencia of Spain,  European Regional Development Fund
under project TEC2005-03575, Catalan Research Council under project 2005SGR00256, and COST action
293 GRAAL, and has been partially done in the context of the {\sc crc Corso} with France Telecom.

\bibliographystyle{abbrv}
\bibliography{MON4}


\begin{appendix}

\section{Small constructions in the proof of Theorem \ref{41}}
\label{ap:1}

\noindent $\mathscr{MON}(3+3,3;4,1)$: $\cB=\{(0,a_0,1;a_2)$, $(1,a_1,2;a_0)$, $(2,a_2,0;a_1)$,
$(a_0,a_1,a_2)\}$.

\medskip
\noindent $\mathscr{MON}(4+4,4;4,1)$: $\cB=\{(1,2,a_3;a_0)$, $(0,3,a_2;a_1)$, $(a_1,1,3;a_0)$,
$(a_0,a_2,1;0)$, $(a_0,a_1,2;0)$, $(a_1,a_3,0;a_0)$, $(2,3,a_3,a_2)\}$.

\medskip
\noindent $\mathscr{MON}(5+5,5;4,1)$: $\cB=\{(1,2,a_3;a_0)$, $(0,3,a_2;a_1)$, $(a_1,1,3;a_0)$,
$(a_0,a_2,1;0)$, $(a_0,a_1,2;0)$, $(a_1,a_3,0;a_0)$, $(2,a_2,4;a_4)$, $(3,a_3,4)$, $(a_2,a_3,a_4)$,
$(2,3,a_4)$, $(0,4,a_0,a_4)$, $(1,4,a_1,a_4)\}$.

\medskip
\noindent $\mathscr{MON}(6+6,6;4,1)$: $\cB=\{(1,2,a_3;a_0)$, $(0,3,a_2;a_1)$, $(a_1,1,3;a_0)$,
$(a_0,a_2,1;0)$, $(a_0,a_1,2;0)$, $(a_1,a_3,0;a_0)$, $(4,5,a_5;a_4)$, $(2,a_2,4;a_4)$,
$(2,3,a_4;5)$, $(3,4,a_3)$, $(a_2,a_3,a_4)$, $(0,4,a_0,a_4)$, $(1,4,a_1,a_4)$, $(0,5,a_0,a_5)$,
$(1,5,a_1,a_5)$, $(2,5,a_2,a_5)$, $(3,5,a_3,a_5)\}$.

\section{Small constructions in the proof of Theorem \ref{mon41small}}
\label{ap:2}

\noindent $\mathscr{MON}(1+2,1;4,1)$:  $\cB=\{(0,a_0,a_1)\}$.

\medskip
\noindent $\mathscr{MON}(2+3,2;4,1)$:  $\cB=\{(0,a_0,a_1)$, $(1,a_1,a_2)$, $(0,1,a_0,a_2) \}$.

\medskip
\noindent $\mathscr{MON}(3+4,3;4,1)$:  $\cB=\{(0,a_0,a_1)$, $(1,a_1,a_2)$, $(0,1,a_0,a_2)$,
$(2,a_2,a_3)$, $(0,2,a_0,a_3)$, $(1,2,a_1,a_3) \}$.

\medskip
\noindent $\mathscr{MON}(4+5,4;4,1)$:  $\cB=\{(0,1,a_0;a_3)$, $(0,2,a_1; a_3)$, $(0,3,a_2;a_3)$,
$(2,3,a_0;a_4)$, $(1,3,a_1;a_4)$, $(1,2,a_3;3)$, $(0,a_3,a_4;3)$, $(1,a_2,a_4;2)$,
$(a_0,a_1,a_2;2)\}$.

\section{Small constructions in the proof of Theorem \ref{42o}}
\label{ap:3}

\noindent $\mathscr{ON}(8,4)$ with 4 triangles: $\cB=\{(1,2,0;4)$, $(0,3,6;7)$, $(0,7,5;2)$,
$(4,5,3;1)$, $(1,4,7)$, $(1,5,6)$, $(2,3,7)$, $(2,4,6)\}$.

\medskip
\noindent $\mathscr{MON}(7+4,7;4,2)$:  $\cB=\{ (a_0,4,2;3)$, $(a_0,3,6;0)$, $(a_0,0,5;1)$,
$(a_1,5,3;4)$, $(a_1,4,6;1)$, $(a_1,1,0;2)$, $(a_2,0,4;5)$, $(a_2,6,5;a_3)$, $(a_2,1,2;5)$,
$(0,3,a_3;2)$, $(1,a_0,a_2,3)$, $(a_0,a_1,a_2,a_3)$, $(a_1,2,6,a_3)$,
 $(1,4,a_3)\}$.

\section{Small constructions in the proof of Theorem \ref{thm:MON43}}
\label{ap:4}

\noindent $\mathscr{MON}(13+3,13;4,3)$: $\cB=\{(5+i,4+i,1+i;a_1)\mid i=0,1,\ldots,9\}\cup
\{(1+i,5+i,4+i;a_0)\mid  i=10,11,12\}\cup \{(3+i,1+i,9+i;a_2)\mid i=6,7,\ldots,12\}\cup
\{(9+i,1+i,3+i;a_0)\mid i=1,2,\ldots,5\}\cup \{(9,3,1;a_2), (0,a_1,a_2;12), (12,a_1,a_0;0),(a_0, 9,
a_2, 10), (a_0,a_2,11;a_1),  (a_0,9,a_2,10)\}$, where the sums are computed modulo 13.

\medskip
\noindent $\mathscr{MON}(15+4,15;4,3)$: $\cB=\{(1,2,3)$, $(a_0,4,a_1,5),$ $(a_0,10,a_1,11)$,
$(5,4,1;a_3)$, $(7,1,6;a_1)$, $(6,4,2;a_3)$, $(7,5,2;a_2)$, $(4,7,3;a_2)$, $(6,5,3;a_3),$
$(9,1,8;a_1)$, $(10,1,14;a_0)$, $(11,1,0;a_2)$, $(13,1,12;a_2)$, $(10,2,8;a_2)$, $(11,2,9;a_0)$,
$(12,2,14;a_2),$ $ (0,2,13;a_3)$, $(8,3,11;a_3)$, $(10,3,12;a_0)$, $(13,3,9;a_2)$, $(14,3,0;a_3)$,
$(12,8,4;a_2)$, $(11,4,13;a_0),$ $ (0,10,4;a_3)$, $(9,4,14;a_1)$, $(8,5,13;a_2)$, $(0,5,12;a_3)$,
$(14,11,5;a_3)$, $(10,9,5;a_2)$, $(8,6,0;a_0),$ $ (14,13,6;a_3)$, $(9,6,12;a_1)$, $(10,6,11;a_2)$,
$(14,8,7;a_3)$, $(9,7,0;a_1)$, $(10,7,13;a_1)$, $(12,11,7;a_0),$ $ (6,a_2,a_0;2)$, $(7,a_2,a_1;3)$,
$(10,a_3,a_2;1)$, $(1,a_0,a_1;2)$, $(9,a_1,a_3;14)$, $(8,a_3,a_0;3)\}$.

\medskip \noindent $\mathscr{MON}(17+3,17;4,3)$: $\cB=\{(7,16,0)$, $(a_0,a_2,0)$, $(a_0,1,2;3)$, $(a_0,3,4;1)$,
$(4,5,2;a_1)$, $(1,3,5;a_0),$ $(16,a_0,a_1;a_2)$, $(6,10,1;a_1)$, $(9,14,1;a_2)$, $(15,1,7;a_2)$,
$(1,8,12;a_2)$, $(1,0,13;a_2)$, $(1,16,11;a_1),$ $(2,11,6;a_1)$, $(2,16,8;a_2)$, $(10,15,2;a_2)$,
$(9,2,13;a_1)$, $(0,2,12;a_1)$, $(2,7,14;a_2)$, $(6,13,3;a_1),$ $(11,3,7;a_1)$, $(12,3,16;a_2)$,
$(9,0,3;a_2)$, $(3,10,14;a_1)$, $(8,3,15;a_1)$, $(14,6,4;a_2)$, $(4,11,15;a_2),$ $(7,12,4;a_1)$,
$(13,4,8;a_1)$, $(4,16,9;a_2)$, $(0,4,10;a_1)$, $(5,12,6;a_2)$, $(7,13,5;a_2)$, $(8,14,5;a_1),$
$(15,5,9;a_1)$, $(5,16,10;a_2)$, $(5,0,11;a_2)$, $(9,7,6;a_0)$, $(10,8,7;a_0)$, $(11,9,8;a_0)$,
$(12,10,9;a_0),$ $(13,11,10;a_0)$, $(14,12,11;a_0)$, $(15,13,12;a_0)$, $(16,14,13;a_0)$,
$(0,15,14;a_0)$, $(6,16,15;a_0),$ $(8,6,0;a_1)\}$.

\medskip
\noindent $\mathscr{MON}(17+4,17;4,3)$: $\cB=\{(2,9,11)$, $(9,12,16),$ $(a_0,13,14;15)$,
$(a_0,15,16;13)$, $(16,0,14;a_1)$, $(13,15,0;a_0)$, $(13,2,1;a_3)$, $(13,12,3;a_3),$
$(13,11,4;a_3)$, $(5,10,13;a_1)$, $(6,9,13;a_2)$, $(7,8,13;a_3)$, $(14,4,2;a_3)$, $(14,12,5;a_3)$,
$(11,14,6;a_3),$ $(14,10,7;a_3)$, $(1,3,14;a_2)$, $(9,8,14;a_3)$, $(1,4,15;a_1)$, $(3,5,15;a_2)$,
$(2,6,15;a_3)$, $(15,7,12;a_3),$ $(15,11,8;a_3)$, $(1,16,5;a_1)$, $(6,4,16;a_2)$, $(3,7,16;a_3)$,
$(2,8,16;a_1)$, $(10,16,11;a_3)$, $(1,6,0;a_1),$ $(4,8,0;a_2)$, $(10,15,9,;a_3)$, $(2,10,0;a_3)$,
$(5,0,7;a_1)$, $(3,0,9;a_1)$, $(12,0,11;a_1)$, $(1,a_0,7;6),$ $(8,6,a_0;a_3)$, $(9,a_0,5;11)$,
$(10,a_0,4;9)$, $(11,a_0,3;10)$, $(2,a_0,12;8)$, $(8,a_1,1;11)$, $(10,a_1,6;3),$ $ (12,4,a_1;a_3)$,
$(3,a_1,2;7)$, $(1,a_2,9;7)$, $(10,a_2,8;3)$, $(11,a_2,7;4)$, $(12,a_2,6;5)$, $(2,a_2,5;8),$ $
(3,a_2,4;5)$, $(a_1,a_0,a_2;a_3)$, $(12,1,10;a_3)\}$.

\end{appendix}
\end{document}